\documentclass[12pt]{article}
\usepackage{amssymb,amsfonts}
\usepackage{amscd}
\usepackage{verbatim}

\newtheorem{theorem}{Theorem}[section]
\newtheorem{lemma}[theorem]{Lemma}
\newtheorem{proposition}[theorem]{Proposition}

\newtheorem{definition}[theorem]{Definition}
\newtheorem{corollary}[theorem]{Corollary}
\newtheorem{remark}[theorem]{Remark}

\newenvironment{proof}{\bf Proof. \rm}{$\Box$}

\newcommand{\be}{\begin{equation}}
\newcommand{\ee}{\end{equation}}

\newcommand{\norm}[1]{\Vert #1 \Vert}

\newcommand{\cN}{\mathcal{N}}

\newcommand{\bZ}{\mathbb{Z}^k}
\newcommand{\bZp}{\mathbb{Z}_+^k}
\newcommand{\bm}{\textbf{m}}
\newcommand{\bp}{\textbf{p}}
\newcommand{\bn}{\textbf{n}}
\newcommand{\bl}{\textbf{l}}
\newcommand{\bq}{\textbf{q}}
\newcommand{\bfe}{\textbf{e}}
\newcommand{\bz}{\textbf{0}}
\newcommand{\bd}{\textbf{d}}

\begin{document}
\title{Regular dilations of representations of product systems}

\author{ Baruch Solel
\\Department of Mathematics\\Technion\\32000 Haifa, Israel
\\e-mail: mabaruch@techunix.technion.ac.il}

\maketitle

\begin{abstract}
We study completely contractive representations of product systems
$X$ of correspondences over the semigroup $\mathbb{Z}_+^k$. We
present a necessary and sufficient condition for such a
representation to have a regular isometric dilation. We discuss
representations that doubly commute and show that these
representations induce completely contractive representations of
the norm closed algebra generated by the image of the Fock
representation of $X$.
\textbf{ key words and phrases.} correspondences, product systems,
 representations, regular dilations, doubly commuting, Nica covariance.

\end{abstract}

\maketitle
\begin{section}{Introduction}\label{intro}
A $C^*$-correspondence $E$ over a $C^*$-algebra $A$ is a (right)
Hilbert $C^*$-module over $A$ that carries also a left action of
$A$ (by adjointable operators). It is also called a Hilbert
bimodule in the literature. A c.c. representation of $E$ on a
Hilbert space $H$ is a pair $(\sigma,T)$ where $\sigma$ is a
representation of $A$ on $H$ and $T:E\rightarrow B(H)$ is a
completely contractive linear map that is also a bimodule map
(that is, $T(a\cdot \xi \cdot b)=\sigma(a)T(\xi)\sigma(b)$ for
$a,b \in A$ and $\xi \in E$). The representation is said to be
isometric (or Toeplitz ) if $T(\xi)^*T(\eta)=\sigma(\langle
\xi,\eta \rangle)$ for every $\xi,\eta \in E$.

 In \cite{P}, Pimsner associated with such a correspondence two
 $C^*$-algebras ($\mathcal{O}(E)$ and $\mathcal{T}(E)$) with certain universal
  properties. In \cite{MS98} we studied the operator
  algebra $\mathcal{T}_+(E)$ (called the tensor algebra)
 which is universal for c.c. representations of $E$.

A product system $X$ of $C^*$-correspondences over a semigroup $P$
is, roughly speaking, a family $\{X_s: s\in P\}$ of
$C^*$-correspondences
(over the same $C^*$-algebra $A$), with $X_e=A$, such that $X_s \otimes X_t$ is
isomorphic to $X_{st}$ for all $s,t \in P\backslash \{e\}$.  (See Section~\ref{prel}
for the precise definition). A c.c. (respectively, isometric)
 representation of $X$ is a family
$\{T_s\}$ such that, for all $s\in P\setminus \{e\}$, $(T_e,T_s)$ is
a c.c. (respectively, isometric) representation of $X_s$  and such that,
 whenever $x\in X_s$ and $y\in X_t$,
$T_{st}(\theta_{s,t}(x\otimes y))=T_s(x)T_t(y)$ (where
$\theta_{s,t}$ is the
isomorphism from $X_s \otimes X_t$ onto $X_{st}$ ).

If $E$ is a $C^*$-correspondence over $A$ then, setting
$X(n)=E^{\otimes n}$ (and $X(0)=A$), we get a product system over
$P=\mathbb{Z}_+$ and every product system over $\mathbb{Z}_+$
arises in this way.

In \cite{F}, Fowler studied product systems over more general
(discrete) semigroups $P$. He proved the existence of a
$C^*$-algebra $\mathcal{T}(X)$ that is universal with respect to
Toeplitz representations. In \cite[Proposition 3.2]{S}, we proved
the existence of an operator algebra $\mathcal{T}_+(X)$ (the
\emph{universal tensor algebra}) which is universal for c.c.
representations of $X$; that is, there is a c.c. representation of
$X$ whose image generates $\mathcal{T}_+(X)$ and every c.c.
representation of $X$ gives rise to a completely contractive
representation of the algebra $\mathcal{T}_+(X)$.

In \cite[Theorem 4.4]{S} we also proved that every c.c.
representation of a product system $X$ over $P=\mathbb{Z}_+^2$ can
be dilated to an isometric representation of $X$. (This was then
used to dilate a pair of commuting CP maps). Specializing to the
case where $A=\mathbb{C}$ and $X(\bn)=\mathbb{C}\;,\bn\in
\mathbb{Z}_+^2$, this result recovers Ando's dilation result (\cite{A}).
 Ando proved that, given a pair  $(T_1,T_2)$ of commuting contractions
in $B(H)$, there is a Hilbert space $K$, containing $H$, and a pair
$(V_1,V_2)$ of commuting isometries in $B(K)$ such that, for all
$\bn=(n_1,n_2)\in \mathbb{Z}_+^2$,
$$ P_HV_1^{n_1}V_2^{n_2}|H=T_1^{n_1}T_2^{n_2}. $$
It is well known (see \cite{Pa} or \cite{Pau}) that such a result
is false, in general, for $\mathbb{Z}_+^k$, $k\geq 3$ (that is,
for $k$-tuples of commuting contractions with $k\geq 3$). Thus, in
particular, Theorem 4.4 of \cite{S}, cannot be proved for product
systems over $\mathbb{Z}_+^k$, $k\geq 3$.

It is known, however, that, if $(T_1,T_2,\ldots,T_k)$ is a
commuting tuple of contractions in $B(H)$ satisfying an additional
condition, then there are isometries $(V_1,V_2,\ldots,V_k)$ (in
$B(K)$ for some Hilbert space $K$ containing $H$) that dilate
$(T_1,T_2,\ldots,T_k)$.(See \cite{B} or \cite[Theorem 9.1]{SzF}).
 The additional condition requires that, for every subset
 $v\subseteq \{1,\ldots,k\}$,
 \be\label{condS}
S(v):=\sum_{u\subseteq v} (-1)^{|u|}(T^{\bfe(u)})^*T^{\bfe(u)}
\geq 0
\ee
where, for $u=\{i_1,\ldots,i_m\}$, $|u|=m$ and
$T^{\bfe(u)}=T_{i_1} \cdots T_{i_m}$. In fact, this condition is a
necessary and sufficient condition to have an isometric dilation
$(V_1,\ldots,V_k)$ with the additional property that, for every
$\bn,\bm \in \bZp$ with $\bn\wedge \bm=\bz$,
$$P_HV^{\bn *}V^{\bm}|H=T^{\bn *}T^{\bm} $$
where $T^{\bn}=\prod T_i^{n_i}$ and $V^{\bn}=\prod V_i^{n_i}$. Such a
dilation is called a \emph{regular dilation}.

In Definition~\ref{reg} we define regular isometric dilations for
c.c. representations of the product system $X$ (over $\bZp$) and,
in Theorem~\ref{regdil}, we prove that a condition similar to
condition (\ref{condS}) is a necessary and sufficient condition
for the existence of an isometric regular dilation. It is also
possible, in this case, to find an isometric regular
representation that is minimal (in an obvious sense) and, in
Proposition~\ref{unique}, we show that such a dilation is unique
up to unitary equivalence.

In the classical case, it is known (\cite[proposition 9.2]{SzF})
that, if the $k$-tuple $(T_1,\ldots,T_k)$ of contractions doubly
commutes (that is, the operators commute and, in addition,
$T_iT_j^*=T_j^*T_i$ for all $i\neq j$), then it satisfies
condition (\ref{condS}) (and, thus, a regular, minimal isometric
dilation exists). It is also known (\cite[Theorem 1]{GS} or
\cite[Theorem 2]{Ti}) that, in this case, the
regular, minimal isometric dilation also doubly commutes.

In Theorem~\ref{doublycomm} we prove a similar result for
representations of $X$. (See Definition~\ref{dc} for the
definition of a doubly commuting representation of a product
system $X$). Then, in Lemma~\ref{dcnc}, we observe that, for an
isometric representation, the doubly commuting condition is
equivalent to a condition known in the literature (e.g. \cite{N},
 \cite{F} or \cite{FR})
as \emph{Nica
covariance}. We then note, using results of \cite{F}, that the
$C^*$-algebra generated by the image of the Fock representation
$L$ on the Fock space $\mathcal{F}(X):=\sum X(\bn)$ is isomorphic
to the algebra $\mathcal{T}_{cov}(X)$. The algebra
$\mathcal{T}_{cov}(X)$ was studied by Fowler in \cite{F} and was
shown there to be universal for Nica-covariant representations
provided $X$ is compactly aligned (Definition~\ref{compal}).
Considering the Banach algebra generated by the image of the Fock
representation $L$ (and writing $\mathcal{T}_{+,c}(X)$ for it), we
use Theorem~\ref{doublycomm} to show, in Corollary~\ref{repn},
that every doubly commuting, c.c. representation $\{T_{\bn}\}$ of $X$ on $H$ gives
rise to a unique completely contractive representation of
$\mathcal{T}_{+,c}(X)$ mapping $L(x)$, for $x\in X(\bn)$, to
$T_{\bn}(x)$. We refer to $\mathcal{T}_{+,c}(X)$ as the
\emph{concrete tensor algebra} associated with $X$.

Recently, $k$-graphs and the $C^*$-algebras associated with them
have been studied extensively. (See \cite{KPa} where these $C^*$-algebras
 were introduced, the survey article \cite{R} and
the references there). Note that every $k$-graph can be defined by
a product system of graphs over $\bZp$ (\cite{RS}). The algebra
$\mathcal{T}_{+,c}$ for such a product system, associated with a
$k$-graph $\Lambda$, is the ``multivariable" analogue of the quiver
algebra of \cite{MS99} and can be referred to as a $k$-quiver
algebra and denoted $\mathcal{T}_{+,c}(\Lambda)$. These algebras
(and their weak closures) were studied in \cite{KP}. In
Subsection~\ref{ex}.4, we discuss the case of a single-vertex
$k$-graph in more details.

The next section is devoted to recalling some preliminary results
and notation. In Section~\ref{rd} we present and prove the main
results of the paper and in Section~\ref{ex} we present some
examples.

\end{section}

\begin{section}{Preliminaries}\label{prel}
We begin by recalling the notion of a $C^*$-correspondence.
For the general theory of Hilbert
$C^*$-modules which we use, we will follow \cite{L}. In
particular, a Hilbert $C^*$-module $E$ over a $C^*$-algebra $A$
will be a \emph{right} Hilbert $C^*$-module. We write
$\mathcal{L}(E)$ for the algebra of continuous, adjointable
$A$-module maps on $E$. It is known to be a $C^*$-algebra.

\begin{definition}\label{corr}
 A $C^*$-correspondence over a $C^*$-algebra $A$ is a
Hilbert $C^*$-module $E$ over $A$ endowed with the structure of a
left $A$-module via a $^*$-homomorphism $\varphi_E :A \rightarrow
\mathcal{L}(E)$.

\end{definition}

When dealing with a specific $C^*$-correspondence $E$ it will be
convenient to write $\varphi$ (instead of $\varphi_E$) or even to
suppress it and write $a\xi$ or $a\cdot \xi$ for $\varphi(a)\xi$.

If $E$ and $F$ are $C^*$-correspondences over $A$, then the
balanced tensor product $E\otimes_A F$ is a $C^*$-correspondence
over $A$. It is defined as the Hausdorff completion of the
algebraic balanced tensor product with the internal inner product
given by
\be\label{tp}
\langle \xi_1 \otimes \eta_1,\xi_2\otimes \eta_2 \rangle=\langle
\eta_1, \varphi_F(\langle \xi_1,\xi_2 \rangle_E)\eta_2 \rangle_F
\ee
for all $\xi_1,\xi_2 \in E$ and $\eta_1,\eta_2 \in F$. The left
and right actions of $a\in M$ are defined by
\be\label{lrp}
\varphi_{E\otimes F}(a)(\xi \otimes \eta)b=\varphi_E(a)\xi \otimes
\eta b
\ee
for all $a,b\in M$, $\xi\in E$ and $\eta\in F$.

\begin{definition}\label{isomcor}
An isomorphism of $C^*$-correspondences $E$ and $F$  is a
surjective, bimodule map that preserves the inner products. We
write $E\cong F$ if such an isomorphism exists.
\end{definition}

If $E$ is a $C^*$-correspondence over $A$ and $\sigma$ is a
representation of $A$ on a Hilbert space $H$  then
$E\otimes_{\sigma}H$ is the Hilbert space obtained as the
Hausdorff completion of the algebraic tensor product with respect
to $\langle \xi \otimes h, \eta \otimes k \rangle=\langle
h,\sigma(\langle \xi,\eta \rangle_E)k\rangle_H $. Given an
operator $X\in \mathcal{L}(E)$ and an operator $S\in \sigma(A)'$,
the map $\xi \otimes h \mapsto X\xi \otimes Sh $ defines a bounded
operator $X\otimes S$ on $E\otimes_{\sigma}H$. When $S=I_E$ and
$X=\varphi_E(a)$ (for $a\in A$) we get a representation of $A$ on
this Hilbert space.
 We frequently
write $a\otimes I_H$ for $\varphi(a)\otimes I_H$.

\begin{definition}
\label{Definition1.12}Let $E$ be a $C^*$-correspondence over a $C^*$-algebra
 $A$. Then a completely contractive covariant
representation of $E$ (or, simply, a c.c. representation of $E$) on a
Hilbert space $H$ is a pair $(\sigma,T)$, where

\begin{enumerate}
\item[(1)] $\sigma$ is a  $\ast$-representation of $A$ in $B(H)$.
\item[(2)] $T$ is a linear, completely contractive map from $E$ to
$B(H)$.
\item[(3)] $T$ is a bimodule map in the sense that $T(a\xi b)=\sigma(a)T(\xi
)\sigma(b)$, $\xi\in E$, and $a,b\in A$.
\end{enumerate}
Such a representation is said to be isometric if, for every
$\xi,\eta\in E$, $T(\xi)^*T(\eta)=\sigma(\langle \xi,\eta
\rangle)$.

\end{definition}

It should be noted that there is a natural way to view $E$ as an operator
space (by viewing it as a subspace of its linking algebra) and this defines
the operator space structure of $E$ to which Definition
\ref{Definition1.12} refers when it is asserted that $T$ is completely contractive.

As we showed in \cite[Lemmas 3.4--3.6]{MS98}, if a completely
contractive covariant representation,
$(\sigma,T)$, of $E$ in $B(H)$ is given, then it determines a contraction
$\tilde{T}:E\otimes_{\sigma}H\rightarrow H$ defined by the formula $\tilde
{T}(\eta\otimes h):=T(\eta)h$, $\eta\otimes h\in E\otimes_{\sigma}H$. The
operator $\tilde{T}$ satisfies
\begin{equation}
\tilde{T}(\varphi(\cdot)\otimes I)=\sigma(\cdot)\tilde{T}.\label{covariance}%
\end{equation}
In fact we have the following lemma from \cite[Lemma 2.16]{MSNP}.

\begin{lemma}
\label{CovRep}The map $(\sigma,T)\rightarrow\tilde{T}$ is a
bijection between all completely contractive covariant
representations $(\sigma,T)$ of $E$ on the Hilbert space $H$ and
contractive operators $\tilde{T}:E\otimes_{\sigma }H\rightarrow H$
that satisfy equation (\ref{covariance}). Given $\sigma$ and a
contraction $\tilde{T}$ satisfying the covariance condition
(\ref{covariance}), we get a  completely contractive covariant
representation $(\sigma,T)$ of $E$ on $H$ by setting $T(\xi)h:=\tilde{T}%
(\xi\otimes h)$.

Moreover, the representation $(\sigma,T)$ is an isometric
representation if and only if $\tilde{T}$ is an isometry.
\end{lemma}

\begin{remark}
\label{GenPowers}In addition to $\tilde{T}$ we also require the
\textquotedblleft generalized higher powers\textquotedblright\emph{\ }of
$\tilde{T}$. These are\emph{\ }maps$\;\tilde{T}_{n}:E^{\otimes n}\otimes
H\rightarrow H\;$defined by the equation$\;\tilde{T}_{n}(\xi_{1}\otimes
\ldots\otimes\xi_{n}\otimes h)=T(\xi_{1})\cdots T(\xi_{n})h$, $\xi_{1}%
\otimes\ldots\otimes\xi_{n}\otimes h\in E^{\otimes n}\otimes H$.
One checks easily that
 $\tilde{T}_n=\tilde{T}\circ
(I_E\otimes \tilde{T})\circ \cdots \circ (I_{E^{\otimes n-1}}
\otimes \tilde{T})$, $n>1$.

\end{remark}

\end{section}
\begin{section}{Regular dilations}\label{rd}
 In the following we follow the
notation of Fowler (\cite{F}).
Let $P$ be the semigroup $\mathbb{Z}_+^k$.
 Suppose $p:X \rightarrow P$
is a family of $C^*$-correspondences over $A$. Write $X(\bn)$ for the
correspondence $p^{-1}(\bn)$ for $\bn=(n_1,\ldots,n_k)\in P$ and
 $\varphi_{\bn} :A
\rightarrow \mathcal{L}(X(\bn))$ for the left action of $A$ on $X(\bn)$.
We say that $X$ is a \emph{product system over} $\bZp$ if $X$ is a
semigroup, $p$ is a semigroup homomorphism and, for each $\bn,\bm \in
\bZp \backslash \{\bz\}$, the map $(x,y) \in X(\bn) \times X(\bm)
\rightarrow
xy \in
X(\bn+\bm)$ extends to an isomorphism $\theta_{\bn,\bm}$ of
correspondences from $X(\bn) \otimes X(\bm)$ onto $X(\bn+\bm)$. We also
require that $X(\bz)=A$ and that the multiplications $X(\bz)\times
X(\bn)
\rightarrow X(\bn)$ and $X(\bn)\times X(\bz) \rightarrow X(\bn)$ are given by
the left and right actions of $A$ on $X(\bn)$.

 The associativity of
the multiplication means that,
for every
$\bn,\bm,\bp \in \bZp$,
\be\label{assoc}
\theta_{\bn+\bm,\bp}(\theta_{\bn,\bm}\otimes
I_{\bp})=\theta_{\bn,\bm+\bp}(I_{\bn}\otimes \theta_{\bm,\bp})
\ee
where, for $\bm\in \bZp$, $I_{\bm}$ stands for the identity of
$X(\bm)$.
We shall write $\bfe_i$ for the element in $\bZp$ whose $i$th
entry is $1$ and all other entries are $0$ and, for a subset
$u\subset \{1,\ldots, k\}$, we write $\bfe(u)=\sum \{\bfe_i: i\in u \}$.

Given a product system $X$ over $\bZp$, we set $E_i=X(\bfe_i)$ for
$1\leq i \leq k$. It will be convenient to write $E_i^n$ for the $n$-fold
 tensor product $E_i^{\otimes n}$ and to identify $X(\bn)$ (for
$\bn\in \bZp$) with $E_1^{n_1}\otimes E_2^{n_2} \otimes \cdots
\otimes E_k^{n_k}$ (where these tensor products are the balanced
tensor products over $A$). That means, in particular, that the
isomorphisms $\theta_{\bfe_i,\bfe_j}$, for $i\leq j$, are identity
maps. Setting $t_{i,j}=\theta_{\bfe_i,\bfe_j}:E_i\otimes E_j \rightarrow
 E_j \otimes E_i$ for $i\geq j$ (and
$t_{i,j}=t_{j,i}^{-1}$ for $i<j$), one can check that the family
$\{t_{i,j}: 1\leq i,j \leq k\}$ satisfies
\be\label{assoct}
(t_{j,i}\otimes
I_{\bfe_l})(I_{\bfe_j}\otimes t_{l,i})(t_{l,j}\otimes I_{\bfe_i})=
(I_{\bfe_i}\otimes t_{l,j})
(t_{l,i}\otimes I_{\bfe_j})(I_{\bfe_l}\otimes t_{j,i})
\ee
for every $1\leq i,j,l \leq k$.
One can also check (but we omit the tedious computation) that,
given $k$ correspondences $E_1,\ldots,E_k$ over the $C^*$-algebra
$A$ and a family $\{t_{i,j}:1\leq i,j \leq k\}$ such that
$t_{i,j}:E_i \otimes E_j \rightarrow E_j\otimes E_i$ is an
isomorphism, $t_{i,j}=t_{j,i}^{-1}$ and $t_{i,i}$ is the identity
map, it determines, in a unique way, a product system $X$ (with
$X(\bn)=E_1^{n_1}\otimes \cdots \otimes E_k^{n_k}$)
whose
isomorphisms $\{\theta_{\bn,\bm}\}$ satisfy
$\theta_{\bfe_i,\bfe_j}=id$ if $i\leq j$ and
$\theta_{\bfe_i,\bfe_j}=
t_{i,j}$ if $i>j$.

\begin{definition}\label{ccrep}
A c.c. representation of $X$ on a Hilbert space $H$ is given by a
non degenerate
representation $\sigma$ of $A$ on $H$ and $k$ completely
contractive maps $T^{(i)}:E_i \rightarrow B(H)$ such that, for each
$1\leq i \leq k$, $(\sigma,T^{(i)})$ is a c.c. representation of $E_i$
and, for $i,j$, they satisfy the commutation relation
\be\label{commute}
\tilde{T}^{(i)}(I_{E_i}\otimes \tilde{T}^{(j)})=\tilde{T}^{(j)}(I_{E_j}\otimes
\tilde{T}^{(i)})\circ (t_{i,j}\otimes I_H)
\ee
\end{definition}

Recall that we write $\tilde{T}^{(i)}_n$ (where $1\leq i \leq k$ and
$n\geq 0$) for $\tilde{T}^{(i)}(I_i\otimes \tilde{T}^{(i)})\cdots
(I_i \otimes I_i \otimes \cdots \otimes \tilde{T}^{(i)}) :E_i^n
\otimes H \rightarrow H$ (where $I_i$ stands for $I_{E_i}$). Similarly,
for $\bn\in \bZp$, we write
\be
\tilde{T}_{\bn}=\tilde{T}^{(1)}_{n_1}(I_{n_1\bfe_1} \otimes
\tilde{T}^{(2)}_{n_2})\cdots (I_{\bn - n_k\bfe_k}\otimes
\tilde{T}^{(k)}_{n_k}):X(\bn)\otimes H \rightarrow H.
\ee
The map $T_{\bn}:X(\bn)\rightarrow B(H)$ is then defined by
$T_{\bn}(\xi)h=\tilde{T}_{\bn}(\xi \otimes h)$ (for $h\in H$).
It follows from (\ref{commute}) that, for $\bn,\bm \in \bZp$, $\xi \in
X(\bn)$
 and $\eta \in X(\bm)$,
\be\label{tmn}
T_{\bn+\bm}(\theta_{\bn,\bm}(\xi \otimes
\eta))=T_{\bn}(\xi)T_{\bm}(\eta).
\ee
 So that Definition~\ref{ccrep}
agrees with the definition stated in Section~\ref{intro}.

For $\bn=(n_1,n_2,\ldots n_k)\in \bZ$ we write $\bn_+$ for the
vector whose $i$th entry is $\max\{n_i,0\}$ and $\bn_{-}$ for
$\bn_+-\bn$. We also write
\be
 T(\bn)=\tilde{T}_{\bn_-}^*\tilde{T}_{\bn_+}:X(\bn_+)\otimes H
\rightarrow X(\bn_-)\otimes H.
\ee
\begin{definition}\label{reg}
Let $(\sigma,T^{(1)}, \ldots, T^{(k)})$ be a c.c. representation
of $X$ on $H$. A regular isometric dilation of $(\sigma,T^{(1)}, \ldots,
T^{(k)})$ is a representation $(\rho,V^{(1)}, \\ \ldots, V^{(k)})$ of
$X$ on a Hilbert space $K$, containing $H$, such that
\begin{enumerate}
\item[(i)] Each $\tilde{V}^{(i)}$ is an isometry (from $E_i
\otimes K$ into $K$).
\item[(ii)] $H$ is invariant for every $V^{(i)}(\xi)^*$, $\xi \in
E_i$.
\item[(iii)] $H$ is reducing for $\rho$ and $\rho(a)|H=\sigma(a)$
for $a\in A$.
\item[(iv)] For every $\bn\in \bZ$, $(I_{X(\bn_-)}\otimes P_H)
V(\bn)|X(\bn_+)\otimes H=T(\bn)$.

\end{enumerate}
Such a dilation is said to be minimal if the smallest closed
subspace of $K$ that contains $H$ and is invariant under all
$V^{(i)}(\xi)$, for $\xi \in E_i$, is $K$.
\end{definition}

Note that the word ``regular" refers to the fact that we require
(iv) to hold for every $\bn\in \bZ$ and not only for $\bn\in
\bZp$.

In the following, in order to avoid cumbersome notation, we shall
often
suppress the isomorphisms between $X(\bn)\otimes X(\bm)$ and
$X(\bn+\bm)$. For example, the map
$I_{\bp-\bfe(u)}\otimes\tilde{T}^*_{\bfe(u)}\tilde{T}_{\bfe(u)}$,
appearing in the statement of Lemma~\ref{comp} below, is a map from
$X(\bp-\bfe(u))\otimes X(\bfe(u))\otimes H$ to itself but we view
it there as a map from $X(\bp)\otimes H$ to itself, invoking these
isomorphisms. Another example is Equation (\ref{tmn}) which will
be frequently used in the form
$$ \tilde{T}_{\bn+\bm}=\tilde{T}_{\bn}(I_{\bn}\otimes
\tilde{T}_{\bm}).$$

The following, technical, lemma will be needed in the proof of the
next theorem.

\begin{lemma}\label{comp}
Let  $(\sigma,\{T^{(i)}\})$ be a c.c. representation of $X$ on $H$.
Write $R=(R(\bp,\bq))_{\bp,\bq \in \bZp}$ for the (infinite,
operator-valued) matrix defined by
$$R(\bp,\bq)=I_{\bq-(\bq-\bp)_+}\otimes T(\bq-\bp) :X(\bq)\otimes H
\rightarrow X(\bp)\otimes H.$$
Write $S=(S(\bp,\bq))_{\bp,\bq \in
\bZp}$ for the matrix defined by $S(\bp,\bq)=R(\bp,\bq)$ if
$\bq\geq \bp$ and $S(\bp,\bq)=0$ otherwise. Also, let $D$ be the
diagonal matrix with
$$D(\bp,\bp)=\sum_{u\subseteq
\{1,\ldots,k\},\bfe(u)\leq \bp}
(-1)^{|u|}(I_{\bp-\bfe(u)}\otimes\tilde{T}^*_{\bfe(u)}\tilde{T}_{\bfe(u)})
:X(\bp)\otimes H \rightarrow X(\bp)\otimes H
$$
 for $\bp\in
\bZp$. Then
\be\label{R}
R=S^*DS.
\ee
Also, let $L$ be the (operator-valued) matrix given by
$L(\bn,\bm)=(-1)^{|v|} I_{\bn}\otimes T(\bfe(v))$ if
$\bm-\bn=\bfe(v)$ and $0$ otherwise. Then $SL=I$ and
\be\label{DL}  D=L^*RL . \ee
\end{lemma}
\begin{remark}\label{infpos}
Before we turn to the proof, note that, although we multiply here
infinite matrices, the sums involved in the computations of the
entries of the product are all finite sums. The precise meaning of
Equation (\ref{R}) is $\langle Rh,g\rangle= \langle DSh,Sg
\rangle$ for $h\in X(\bp)\otimes H$ and $g\in X(\bq)\otimes H$.
 Thus, it holds for all
  $h,g$ in the vector space $\mathcal{H}_0$, which is the (algebraic) sum
   $\sum_{\bp\in \bZp} X(p)\otimes H$.
A similar remark applies to Equation (\ref{DL}). It
thus follows from the lemma that, $R$ is positive on this space
(in the sense that $\langle Rh,h\rangle\geq 0$ for every $h\in \mathcal{H}_0$) if and
only if $D$ is positive (in a similar sense).
\end{remark}

\begin{proof} (Of Lemma~\ref{comp})
Given $\bz\neq \bn \in \bZp$, it is easy to check that
\be\label{sum}
\sum_{u\subseteq \{1,\ldots,k\}, \bfe(u)\leq \bn} (-1)^{|u|} =0.
\ee
If $\bn=\bz$, this sum is, of course, $1$.

Now compute, for $\bp,\bq \in \bZp$,
$$(S^*DS)(\bp,\bq)=\sum_{\bl\leq \bp\wedge \bq}
S(\bl,\bp)^*D(\bl,\bl)S(\bl,\bq)=$$
$$\sum_{\bfe(u)\leq \bl\leq \bp\wedge \bq}(-1)^{|u|}(I_{\bl} \otimes T(\bp -\bl)^*)
(I_{\bl-\bfe(u)}\otimes\tilde{T}^*_{\bfe(u)}\tilde{T}_{\bfe(u)})(I_{\bl}\otimes
T(\bq-\bl))=$$
$$\sum_{\bfe(u)\leq \bl\leq \bp\wedge \bq}(-1)^{|u|}(I_{\bl} \otimes
\tilde{T}_{\bp -\bl}^*)
(I_{\bl-\bfe(u)}\otimes\tilde{T}^*_{\bfe(u)}\tilde{T}_{\bfe(u)})(I_{\bl}\otimes
\tilde{T}_{\bq-\bl})=$$
$$\sum_{\bfe(u)\leq \bl\leq \bp\wedge
\bq}(-1)^{|u|}I_{\bl-\bfe(u)}\otimes
(\tilde{T}_{\bfe(u)+\bp-\bl}^*\tilde{T}_{\bfe(u)+\bq-\bl})=$$
$$\sum_{\bz\leq \bm \leq \bp\wedge \bq} \sum_{\bfe(u)+\bm\leq
\bp\wedge \bq} (-1)^{|u|}(I_{\bm}\otimes
\tilde{T}_{\bp+\bm}^*\tilde{T}_{\bp+\bm})=$$
$$\sum_{\bz\leq \bm \leq \bp\wedge \bq} (\sum_{\bfe(u)+\bm\leq
\bp\wedge \bq} (-1)^{|u|})(I_{\bm}\otimes
\tilde{T}_{\bp+\bm}^*\tilde{T}_{\bp+\bm}).$$
Applying (\ref{sum}), the last sum is equal to
$$I_{\bp\wedge \bq}\otimes (\tilde{T}_{\bp-(\bp\wedge \bq)}^*\tilde{T}_{\bq-(\bp\wedge \bq)})
=I_{\bq-(\bq-\bp)_+}\otimes T(\bq-\bp).$$
This proves that $R=S^*DS$.

Now, let $L$ be as in the statement of the lemma and compute
$$(SL)(\bp,\bq)=\sum_{\bp \leq \bl \leq \bq}
S(\bp,\bl)L(\bl,\bq)=$$
$$\sum_{\bp\leq \bl,\bq=\bl+\bfe(v)}
(-1)^{|v|}(I_{\bp}\otimes T(\bl-\bp))(I_{\bl}\otimes
T(\bfe(v)))=$$
$$\sum_{\bp\leq
\bl,\bq=\bl+\bfe(v)}(-1)^{|v|}I_{\bp}\otimes
(T(\bl-\bp)(I_{\bl-\bp}\otimes T(\bfe(v))))=\sum (-1)^{|v|}(I_{\bp}\otimes
T(\bq-\bp)) $$
where the last sum runs over all $v\subseteq \{1\leq i \leq k:
p_i<q_i \}$. The argument at the beginning of the proof shows that
this is non zero only if $\bp=\bq$ and, in that case, it is equal
to $I_{\bq}$. This shows that $SL=I$ and, consequently, $D=L^*RL$.

\end{proof}

\begin{theorem}\label{regdil}
A c.c. representation $(\sigma,\{T^{(i)}\})$ of $X$ on $H$ has a
regular isometric dilation if and only if, for every
$v\subseteq \{1,\ldots,k \}$,
\be\label{D}
\sum_{u\subseteq v}
(-1)^{|u|}(I_{\bfe(v)-\bfe(u)}\otimes \tilde{T}^*_{\bfe(u)}\tilde{T}_{\bfe(u)}) \geq 0
\ee
where $|u|$ is the number of elements in $u$.

The regular isometric dilation, when it exists, can be chosen minimal.
\end{theorem}
\begin{proof}
Suppose condition (\ref{D}) holds.
Write $\mathcal{H}_0$ for the vector space of all finitely
supported functions $g$ on $\bZp$ with $g(\bm)\in X(\bm)\otimes H$
for all $\bm \in \bZp$. On $\mathcal{H}_0$ we consider the
following sesquilinear form
\be\label{sesq}
\langle g,f \rangle= \sum_{\bn,\bm \geq \bz} \langle
(I_{\bn-(\bn-\bm)_+}\otimes T(\bn-\bm))g(\bn),f(\bm)\rangle.
\ee
Lemma~\ref{comp} (together with condition (\ref{D})) implies that
this form is positive semidefinite.
Let $\cN$ be the space of all $g\in \mathcal{H}_0$ with $\langle
g,g \rangle=0$ and write $K$ for the Hilbert space obtained by
completing the quotient $\mathcal{H}_0/ \cN$ with respect to the
inner product defined by (\ref{sesq}).

We first embedd $H$ into $K$. For that, define $W:H\rightarrow K$
by $Wh=h\delta_{\bz} +\cN$ where $h\delta_{\bz}(\bz)=h\in
H=X(\bz)\otimes H$ and $h\delta_{\bz}(\bn)=0$ if $\bn\neq \bz$.
Then, for $h,f \in H$, $\langle Wh,Wf \rangle=\langle
h\delta_{\bz},f\delta_{\bz}\rangle= \langle h,f \rangle$. Thus $W$
is an isometry of $H$ into $K$.

Now, for $a\in A$ and $g \in \mathcal{H}_0$, we set
$\rho(a)(g+\cN)=f+\cN$ where $f(\bm)=(\varphi_{X(\bm)}(a)\otimes
I_H)g(\bm)$.
Note that, if $a\in A$ and $\bn,\bm\in \bZp$ satisfy $\bn \wedge
\bm \neq \bz$, then $(I_{\bn \wedge \bm}\otimes
T(\bn-\bm))(\varphi_{X(\bn)}(a)\otimes I_H)=(\varphi_{X(\bm)}(a)\otimes I_H)
(I_{\bn \wedge \bm}\otimes
T(\bn-\bm))$ since $\varphi(a)$ acts on the left most factor in
$X(\bn \wedge \bm)$. If $\bn\wedge \bm=\bz$ we still have the same
equality since, in this case, $T(\bn-\bm)(\varphi_{X(\bn)}(a)\otimes
I_H)=\tilde{T}^*_{\bm}\tilde{T}_{\bn}(\varphi_{X(\bn)}(a)\otimes
I_H)=\tilde{T}^*_{\bm}\sigma(a)\tilde{T}_{\bn}=(\varphi_{X(\bm)}(a)\otimes
I_H)\tilde{T}^*_{\bm}\tilde{T}_{\bn}$. Thus, letting $C(a)$ be the
diagonal matrix with $\varphi_{X(\bn)}\otimes I_H$ in the
$\bn,\bn$ entry, we find that $C(a)$ commutes with $R$ (where $R$
is as in Lemma~\ref{comp}). Clearly $\norm{C(a)}\leq \norm{a}$
and, therefore, $C(a)^*RC(a)\leq \norm{a}^2 R$. It follows that
the map $\rho(a)$, defined above, is a well defined bounded
operator on $K$. It is easy to check that $\rho$ is indeed a
$C^*$-representation of $A$ on $K$.

For $1\leq i \leq k$ and $\xi \in E_i$, define $V^{(i)}(\xi)(g+\cN)
=g_i+\cN$ where $g_i(\bn)=\xi \otimes g(\bn - \bfe_i)$ if $\bn\geq
\bfe_i$ and is $0$ otherwise.

Fix $a,b \in A$ and $\xi \in E_i$ and write $f+ \cN$ for
$V^{(i)}(\varphi_{E_i}(a)\xi b)(g+\cN)$. Then
$f(\bn)=\varphi_{E_i}(a)\xi b \otimes
g(\bn-\bfe_i)=(\varphi_{X(\bn)}(a)\otimes I_H)(\xi b\otimes
g(\bn-\bfe_i))=(\varphi_{X(\bn)}(a)\otimes I_H)(\xi \otimes
(\varphi_{X(\bn-\bfe_i)}(b)\otimes I_H)
g(\bn-\bfe_i))=\rho(a)(\xi \otimes (\rho(b)g)(\bn-\bfe_i))$. Thus
$V^{(i)}(\varphi_{E_i}(a)\xi
b)(g+\cN)=\rho(a)V^{(i)}(\xi)\rho(b)(g+\cN)$ proving the
covariance property of $V^{(i)}$. We now turn to show that
$\tilde{V}^{(i)}$ is an isometry.

For this, fix $\xi,\eta \in E_i$ and $g,f \in \mathcal{H}_0$,
write $V^{(i)}(\xi)(g+\cN)=g_i+\cN$, $V^{(i)}(\xi)(f+\cN)=f_i+\cN$
and compute
$$\langle g_i,f_i\rangle = \sum_{\bn,\bm \geq \bfe_i}\langle(I_{\bn\wedge
\bm} \otimes T(\bn-\bm))(\xi \otimes g(\bn-\bfe_i)),\eta \otimes
f(\bm-\bfe_i)\rangle =$$
$$\sum_{\bn,\bm \geq \bfe_i} \langle \xi \otimes (I_{\bn\wedge \bm
-\bfe_i}\otimes T(\bn-\bm))g(\bn-\bfe_i),\eta \otimes f(\bm
-\bfe_i) \rangle =$$
$$\sum_{\bn,\bm \geq \bfe_i}\langle  (I_{\bn\wedge \bm
-\bfe_i}\otimes
T(\bn-\bm))g(\bn-\bfe_i),(\varphi_{X(\bm-\bfe_i)}(\langle \xi,\eta
\rangle )\otimes I_H)f(\bm-\bfe_i) \rangle =$$
$$\langle g,\rho(\langle
\xi,\eta \rangle)f \rangle.$$
Thus $V^{(i)}(\xi)^*V^{(i)}(\eta)=\rho(\langle \xi,\eta \rangle)$
so that, for each $1\leq i \leq k$, $(\rho,V^{(i)})$ is an
isometric representation of $E_i$.

Now, for $g\in \mathcal{H}_0$, $\xi \in E_i$ and $h\in H$, we
compute
$$ \langle g, V^{(i)}(\xi)^*Wh \rangle = \langle V^{(i)}(\xi)g,Wh
\rangle = \sum_{\bn\geq \bfe_i}\langle \tilde{T}_{\bn}(\xi \otimes
g(\bn - \bfe_i)),h \rangle =$$
$$ \sum_{\bn\geq \bfe_i}\langle
T^{(i)}(\xi)\tilde{T}_{\bn-\bfe_i}(g(\bn-\bfe_i)),h \rangle =
 \sum_{\bn\geq \bfe_i}\langle
\tilde{T}_{\bn-\bfe_i}(g(\bn-\bfe_i)),T^{(i)}(\xi)^*h \rangle =$$
$$\langle g, WT^{(i)}(\xi)^*h \rangle .$$
Thus $V^{(i)}(\xi)^*W= WT^{(i)}(\xi)^*$. This proves property (ii)
of Definition~\ref{reg}. Property (iii) is easy to check and we
need only to verify (iv).

Note first that, for $\bp \in \bZp$, $\xi \in X(\bp)$ and $g\in
\mathcal{H}_0$, it follows from the definition of $V^{(i)}$ above
that $\tilde{V}_{\bp}(\xi \otimes g)(\bn)=\xi \otimes g(\bn -
\bp)$ if $\bn \geq \bp$ (and it is equal to $0$ otherwise). Thus,
for $h\in H$, $\tilde{V}_{\bp}(\xi \otimes Wh)(\bn)=\xi \otimes h$
if $\bn=\bp$ and $0$ otherwise. Therefore, for $\bn\in \bZ$,
$\xi\in X(\bn_+)$, $\eta\in X(\bn_-)$ and $h_1,h_2 \in H$,
$$\langle V(\bn)(\xi \otimes Wh_1),\eta \otimes Wh_2 \rangle
=\langle \tilde{V}_{\bn_+}(\xi \otimes
Wh_1),\tilde{V}_{\bn_-}(\eta \otimes Wh_2)\rangle =$$
$$ \langle
T(\bn)(\xi \otimes h_1),\eta \otimes h_2 \rangle.$$
This proves that this is indeed a regular isometric dilation.

Now assume that $(\sigma,\{T^{(i)}\})$ has an isometric regular
dilation $(\rho,\{V^{(i)}\})$ (on $K$). Let $R_V,S_V$ and $D_V$ be the
matrices described in Lemma~\ref{comp} with $V$ replacing $T$.
Since $(\rho,\{V^{(i)}\})$ is an isometric representation, it follows
that, for $u\subseteq \{1,\ldots,k\}$,
$\tilde{V}^*_{\bfe(u)}\tilde{V}_{\bfe(u)}$ is the identity map on
$X(\bfe(u))\otimes H$. The argument in the first paragraph of the
proof of Lemma~\ref{comp} now shows that $D_V$ is the identity
matrix and, thus, $R_V=S_V^*D_VS_V \geq 0$. But, since the
dilation is regular, the matrix $R$ (as in Lemma~\ref{comp}) is a
compression of $R_V$. It follows that $R\geq 0$ and, using
Lemma~\ref{comp} again, $D\geq 0$. From this, (\ref{D}) follows.

If a regular, isometric, dilation exists, we can restrict it to
the minimal closed subspace containing $H$ and invariant under all
$V^{(i)}(\xi),\;\xi\in E_i$, to get a minimal one.
\end{proof}

The following lemma is easy to verify but will be useful.

\begin{lemma}\label{isom}
If $(\sigma,\{V^{(i)}\})$ is an isometric representation (that is,
each $\tilde{V}^{(i)}$ is an isometry), then, for $\bn,\bm \in
\bZp$, $\tilde{V}_{\bm}^*\tilde{V}_{\bn}=I_{\bn\wedge \bm}\otimes
V(\bn-\bm)$.
\end{lemma}
\begin{proof}
Compute $\tilde{V}_{\bm}^*\tilde{V}_{\bn}=(I_{\bn\wedge \bm}\otimes
\tilde{V}_{\bm-\bm\wedge \bn}^*)\tilde{V}_{\bm\wedge
\bn}^*\tilde{V}_{\bm\wedge \bn}(I_{\bn\wedge \bm}\otimes
\tilde{V}_{\bn-\bm\wedge \bn})=I_{\bm\wedge \bn}\otimes
V(\bn-\bm)$.
\end{proof}

\begin{proposition}\label{unique}
A minimal, regular, isometric dilation of $(\sigma,\{T^{(i)}\})$ is
unique up to unitary equivalence.
\end{proposition}
\begin{proof}
Suppose $(\rho,\{V^{(i)}\})$ and $(\tau,\{U^{(i)}\})$ are minimal regular
isometric dilations of $(\sigma,\{T^{(i)}\})$ on $K$ and $G$
respectively.
For every $\bn \in \bZp$ write
$K(\bn)=\tilde{V}_{\bn}(X(\bn)\otimes H)$ and $G(\bn)=\tilde{U}_{\bn}(X(\bn)\otimes
H)$ (and, for $\bn=\bz$, $K(\bz)=H=G(\bz)$).
Now, let $R(\bn):K(\bn) \rightarrow G(\bn)$ be defined by
$R(\bn)\tilde{V}_{\bn}(\xi \otimes h)=\tilde{U}_{\bn}(\xi \otimes
h)$ (for $\xi \in X(\bn)$ and $h\in H$) and $R(\bz)=I_H$. For $\bn,\bm \in \bZp$,
$\xi\in X(\bn)$, $\eta\in X(\bm)$ and $h,g \in H$, we have

$$\langle \tilde{V}_{\bn}(\xi\otimes h),\tilde{V}_{\bm}(\eta
\otimes g)\rangle=\langle \tilde{V}_{\bm}^*\tilde{V}_{\bn}(\xi\otimes h),\eta
\otimes g\rangle=$$
$$\langle (I_{\bm\wedge \bn}\otimes V(\bn-\bm))(\xi \otimes
h),\eta \otimes g \rangle=\langle (I_{\bm\wedge \bn}\otimes T(\bn-\bm))(\xi \otimes
h),\eta \otimes g \rangle$$
 where the second equality follows from Lemma~\ref{isom} and
 last one follows from
Definition~\ref{reg} (iv). A similar computation holds for $U$, in
place of $V$, and we get $\langle
R(\bn)k_{\bn},R(\bm)k_{\bm}\rangle =\langle
k_{\bn},k_{\bm}\rangle$ for every $k_{\bn}\in K(\bn)$ and
$k_{\bm}\in K(\bm)$. This shows that each $R(\bn)$ is well defined
and isometric and, also, that there is a unitary operator
 $R:K\rightarrow G$ such that $R|K(\bn)=W(\bn)$ for $\bn \in
 \bZp$.
Fix $1\leq i \leq k$, $\bn\in \bZp$, $\eta \in X(\bn)$, $\xi
\in E_i$ and $h \in H$. Then
$$RV^{(i)}(\xi)\tilde{V}_{\bn}(\eta \otimes
h)=R\tilde{V}_{\bfe_i}(I_{\bfe_i}\otimes \tilde{V}_{\bn})(\xi
\otimes \eta \otimes h)=
R\tilde{V}_{\bn+\bfe_i}(\xi \otimes \eta \otimes h)=$$
$$\tilde{U}_{\bn+\bfe_i}(\xi \otimes \eta \otimes h)
=\tilde{U}_{\bfe_i}(I_{\bfe_i}\otimes \tilde{U}_{\bn})(\xi
\otimes \eta \otimes h)=U^{(i)}(\xi)\tilde{U}_{\bn}(\eta \otimes
h).$$
  It follows from the minimality assumption that, for all $1\leq
 i \leq k$ and $\xi \in E_i$, $RV^{(i)}(\xi)=U^{(i)}(\xi)R$. Similarly, one
 checks that,
 for $a\in A$,
 $R\rho(a)=\tau(a)R$.
\end{proof}

\begin{definition}\label{dc}
We say that a representation $(\sigma,\{T^{(i)}\})$ is a doubly
commuting representation if, for every $i\neq j$ (in
$\{1,\ldots,k\}$), we have
\be\label{eqdc}
\tilde{T}^{(j)*}\tilde{T}^{(i)}=(I_{\bfe_j}\otimes
\tilde{T}^{(i)})(I_{\bfe_i}\otimes \tilde{T}^{(j)*}).
\ee
More precisely, $\tilde{T}^{(j)*}\tilde{T}^{(i)}=(I_{\bfe_j}\otimes
\tilde{T}^{(i)})(t_{i,j}\otimes I_H)(I_{\bfe_i}\otimes
\tilde{T}^{(j)*})$ where $t_{i,j}:E_i\otimes E_j \rightarrow E_j
\otimes E_i$ is the isomorphism as in Equation (\ref{assoct}).
\end{definition}

\begin{lemma}\label{consdc}
Let $(\sigma,T^{(i)})$ be a doubly commuting representation. Then
\begin{enumerate}
\item[(i)] For $\bn,\bm \in \bZp$ with $\bn \wedge \bm=\bz$,
$$(I_{\bm}\otimes \tilde{T}_{\bn})(I_{\bn}\otimes
\tilde{T}_{\bm}^*)=\tilde{T}_{\bm}^*\tilde{T}_{\bn}.$$
In particular, for $\bp\in \bZ$,
$$(I_{\bp_-}\otimes \tilde{T}_{\bp_+})(I_{\bp_+}\otimes
\tilde{T}_{\bp_-}^*)=T(\bp).$$
\item[(ii)] If $\bp,\bq,\bn \in \bZp$ with $\bp \leq \bn$ and
$\bq\wedge \bp =\bz$, then
$$(I_{\bn-\bp+\bq}\otimes \tilde{T}_{\bp}^*\tilde{T}_{\bp})
(I_{\bn}\otimes
\tilde{T}_{\bq}^*\tilde{T}_{\bq})=I_{\bn-\bp}\otimes
\tilde{T}_{\bp+\bq}^*\tilde{T}_{\bp+\bq}.$$
\item[(iii)] For $u\subseteq v \subseteq \{1,\ldots, k\}$ and
$l\notin v$,
$$(I_{\bfe(v)-\bfe(u)+\bfe_l}\otimes
\tilde{T}_{\bfe(u)}^*\tilde{T}_{\bfe(u)})(I_{\bfe(v)+\bfe_l}\otimes
I_H -(I_{\bfe(v)}\otimes \tilde{T}^{(l)*}\tilde{T}^{(l)}))=$$
$$(I_{\bfe(v)-\bfe(u)+\bfe_l}\otimes
\tilde{T}_{\bfe(u)}^*\tilde{T}_{\bfe(u)})-(I_{\bfe(v)-\bfe(u)}\otimes
\tilde{T}^*_{\bfe(u)+\bfe_l}\tilde{T}_{\bfe(u)+\bfe_l}).$$
\item[(iv)] Let $j\neq l$ in $\{1,\ldots,k\}$ and
$\{j,l\}\subseteq w \subseteq \{1,\ldots, k\}$. Then
$$ (I_{\bfe(w)-\bfe_j}\otimes \tilde{T}^{(j)*}\tilde{T}^{(j)})
(I_{\bfe(w)-\bfe_l}\otimes \tilde{T}^{(l)*}\tilde{T}^{(l)})=
I_{\bfe(w)-\bfe_l-\bfe_j}\otimes
\tilde{T}^{*}_{\bfe_j+\bfe_l}\tilde{T}_{\bfe_j+\bfe_l}$$
$$=
(I_{\bfe(w)-\bfe_l}\otimes \tilde{T}^{(l)*}\tilde{T}^{(l)})
(I_{\bfe(w)-\bfe_j}\otimes \tilde{T}^{(j)*}\tilde{T}^{(j)}).$$

\end{enumerate}
\end{lemma}
\begin{proof}
We start by proving part (i) for all $(\bn,\bm)$ with $\bn \wedge
\bm =\bz$ by induction on $r(\bn,\bm):=\sum_i n_i +\sum_j m_j$. If
$r(\bn,\bm)\leq 2$, then either $\bm=\bfe_j$ and $\bn=\bfe_i$ (with $i\neq
j$) and, in this case, (i) follows from the definition, or one of the
tuples is $\bz$ and, in that case, (i) is trivial. Now
assume $\bn,\bm \in \bZp$ with $\bn\wedge \bm=\bz$ and (i) holds for
all $\bp,\bq \in \bZp$ with $\bp\wedge \bq =\bz$ and
$r(\bp,\bq)<r(\bn,\bm)$.
Since now $r(\bn,\bm)>2$, we can find either some $j\in
\{1,\ldots,k\}$ such that $\bfe_j\leq \bn$ and $\bfe_j\neq \bn$ or
some $i$ such that $\bfe_i\leq \bm$ and $\bfe_i\neq \bm$. Assume,
without loss of generality that, for some $j$, $\bfe_j\lvertneqq
\bn$.  Then
$\tilde{T}_{\bn}=\tilde{T}_{\bfe_j}(I_{\bfe_j}\otimes
\tilde{T}_{\bn-\bfe_j})$ and
$$(I_{\bm}\otimes \tilde{T}_{\bn})(I_{\bn}\otimes
\tilde{T}_{\bm}^*)=(I_{\bm}\otimes \tilde{T}_{\bfe_j})
(I_{\bm+\bfe_j}\otimes \tilde{T}_{\bn-\bfe_j})(I_{\bn}\otimes
\tilde{T}_{\bm}^*)=$$
$$(I_{\bm}\otimes
\tilde{T}_{\bfe_j})(I_{\bfe_j}\otimes (I_{\bm}\otimes
\tilde{T}_{\bn-\bfe_j})(I_{\bn-\bfe_j}\otimes
\tilde{T}^*_{\bm})).$$
Using the induction hypothesis, this is equal to
$$(I_{\bm}\otimes
\tilde{T}_{\bfe_j})(I_{\bfe_j}\otimes
\tilde{T}^*_{\bm}\tilde{T}_{\bn-\bfe_j})=(I_{\bm}\otimes
\tilde{T}_{\bfe_j})(I_{\bfe_j}\otimes
\tilde{T}^*_{\bm})(I_{\bfe_j}\otimes
\tilde{T}_{\bn-\bfe_j}).$$ Using the induction hypothesis again
(for $(\bm,\bfe_j)$), we see that this is equal to
$$\tilde{T}^*_{\bm}\tilde{T}_{\bfe_j}(I_{\bfe_j}\otimes
\tilde{T}_{\bn-\bfe_j})=\tilde{T}^*_{\bm}\tilde{T}_{\bn}.$$ This
completes the proof of (i).

To prove (ii), we compute
$$(I_{\bn-\bp+\bq}\otimes \tilde{T}_{\bp}^*\tilde{T}_{\bp})
(I_{\bn}\otimes
\tilde{T}_{\bq}^*\tilde{T}_{\bq})=$$
$$(I_{\bn-\bp+\bq}\otimes \tilde{T}_{\bp}^*)
(I_{\bn-\bp+\bq}\otimes \tilde{T}_{\bp})
(I_{\bn}\otimes
\tilde{T}_{\bq}^*)(I_{\bn}\otimes
\tilde{T}_{\bq})=$$
$$(I_{\bn-\bp+\bq}\otimes \tilde{T}_{\bp}^*)
(I_{\bn-\bp}\otimes (I_{\bq}\otimes
\tilde{T}_{\bq})(I_{\bp}\otimes \tilde{T}^*_{\bq}))
(I_{\bn}\otimes
\tilde{T}_{\bq}).$$ Using part (i), this is equal to
$$(I_{\bn-\bp+\bq}\otimes \tilde{T}_{\bp}^*)
(I_{\bn-\bp}\otimes \tilde{T}_{\bq}^*\tilde{T}_{\bp})
(I_{\bn}\otimes
\tilde{T}_{\bq})=$$
$$(I_{\bn-\bp+\bq}\otimes \tilde{T}_{\bp}^*)
(I_{\bn-\bp}\otimes \tilde{T}^*_{\bq})(I_{\bn-\bp}\otimes
\tilde{T}_{\bp})
(I_{\bn}\otimes
\tilde{T}_{\bq})=$$
$$(I_{\bn-\bp}\otimes \tilde{T}_{\bp+\bq}^*)
(I_{\bn-\bp}\otimes \tilde{T}_{\bp+\bq})=
(I_{\bn-\bp}\otimes \tilde{T}_{\bp+\bq}^*\tilde{T}_{\bp+\bq})$$
completing the proof of (ii).
To prove (iii), apply (ii) with $\bp=\bfe(u)$, $\bq=\bfe_l$ and
$\bn=\bfe(v)$. Part (iv) is also a consequence of (ii). Simply set
$\bp=\bfe_j$, $\bq=\bfe_i$ and $\bn=\bfe(w)-\bfe_l$ to get one
equality and exchange $j$ and $l$ to get the other one.

\end{proof}

\begin{theorem}\label{doublycomm}
If the representation $(\sigma,\{T^{(i)}\})$ is doubly commuting then
it has a regular isometric dilation and the regular isometric
dilation that is minimal is doubly commuting.
\end{theorem}
\begin{proof}
To show that it has a regular isometric dilation, we should verify
condition (\ref{D}) of Theorem~\ref{regdil}.
In fact, we claim that, for every $v\subseteq \{1,\ldots,k\}$, we
have
\be\label{condD}
\sum_{u\subseteq v}
(-1)^{|u|}(I_{\bfe(v)-\bfe(u)}\otimes \tilde{T}^*_{\bfe(u)}\tilde{T}_{\bfe(u)})
=\prod_{i\in v} (I_{\bfe(v)}\otimes I_H-(I_{\bfe(v)-\bfe_i}\otimes
 \tilde{T}^{(i)*}\tilde{T}^{(i)})).
\ee
Since, by Lemma~\ref{consdc} (iv), the operators in the product commute,
this will show that
the condition of Theorem~\ref{regdil} holds.

 We shall prove the claim by
induction on the number of elements in $v$. If $|v|=2$, we can
write $v=\{j,l\}$ and then the claim follows easily from
Lemma~\ref{consdc} (iv). Now assume we know it for $v$ and
$w=v\cup \{l\}$ where $l\notin v$. Tensoring (\ref{condD}) (for $v$) by
$I_{\bfe_l}$, we get
$$\sum_{u\subseteq v}
(-1)^{|u|}(I_{\bfe(w)-\bfe(u)}\otimes \tilde{T}^*_{\bfe(u)}\tilde{T}_{\bfe(u)})
=\prod_{i\in v} (I_{\bfe(w)}\otimes I_H-(I_{\bfe(w)-\bfe_i}\otimes
 \tilde{T}^{(i)*}\tilde{T}^{(i)})).$$ Thus
$$\prod_{i\in w} (I_{\bfe(w)}\otimes I_H-(I_{\bfe(w)-\bfe_i}\otimes
 \tilde{T}^{(i)*}\tilde{T}^{(i)}))=$$
 $$(\sum_{u\subseteq v}
(-1)^{|u|}(I_{\bfe(w)-\bfe(u)}\otimes
\tilde{T}^*_{\bfe(u)}\tilde{T}_{\bfe(u)}))(I_{\bfe(w)}\otimes I_H-(I_{\bfe(w)-\bfe_l}\otimes
 \tilde{T}^{(l)*}\tilde{T}^{(l)})).$$
Using Lemma~\ref{consdc} (iii), this is equal to
$$\sum_{u\subseteq v}(-1)^{|u|}((I_{\bfe(v)-\bfe(u)+\bfe_l}\otimes
\tilde{T}_{\bfe(u)}^*\tilde{T}_{\bfe(u)})-(I_{\bfe(v)-\bfe(u)}\otimes
\tilde{T}^*_{\bfe(u)+\bfe_l}\tilde{T}_{\bfe(u)+\bfe_l}))=$$
$$\sum_{u\subseteq v}(-1)^{|u|}((I_{\bfe(w)-\bfe(u)}\otimes
\tilde{T}_{\bfe(u)}^*\tilde{T}_{\bfe(u)})-(I_{\bfe(w)-\bfe_l-\bfe(u)}\otimes
\tilde{T}^*_{\bfe(u)+\bfe_l}\tilde{T}_{\bfe(u)+\bfe_l}))=$$
$$\sum_{u\subseteq w}(-1)^{|u|}(I_{\bfe(w)-\bfe(u)}\otimes
\tilde{T}^*_{\bfe(u)}\tilde{T}_{\bfe(u)}).$$
This completes the proof of the claim and shows that the
representation has an isometric regular dilation. In this case, it
has an isometric regular dilation $(\rho,\{V^{(i)}\})$ (on $K$) that is minimal
 in the sense that
\be
\bigvee \{\tilde{V}_{\bn}(X(\bn)\otimes H) : \;\bn \in \bZp \}=K.
\ee
To prove that the representation $(\rho,\{V^{(i)}\})$
is doubly commuting, we fix
$i\neq j$ and we should prove the equality
$$\tilde{V}^{(j)*}\tilde{V}^{(i)}=(I_{\bfe_j}\otimes
\tilde{V}^{(i)})(I_{\bfe_i}\otimes \tilde{V}^{(j)*}).$$
On both sides of this equality we have operators from $E_i \otimes
K$ to $E_j\otimes K$. It follows from the minimality condition
that
\be\label{min}
\bigvee \{(I_{\bfe_i}\otimes \tilde{V}_{\bn})(X(\bn+\bfe_i)\otimes H) :
 \;\bn \in \bZp \}=E_i \otimes K.
 \ee
Thus, it suffices to show that, for every $\bn,\bm \in \bZp$, $\xi
\in X(\bn+\bfe_i)$, $\eta \in X(\bm+\bfe_j)$ and $h,g \in H$,
\be
\langle \tilde{V}^{(j)*}\tilde{V}^{(i)}(I_{\bfe_i}\otimes
\tilde{V}_{\bn})(\xi\otimes h),(I_{\bfe_j}\otimes
\tilde{V}_{\bm})(\eta\otimes g) \rangle=
\ee
$$\langle (I_{\bfe_j}\otimes
\tilde{V}^{(i)})(I_{\bfe_i}\otimes \tilde{V}^{(j)*})
(I_{\bfe_i}\otimes
\tilde{V}_{\bn})(\xi\otimes h),(I_{\bfe_j}\otimes
\tilde{V}_{\bm})(\eta\otimes g) \rangle.$$
The left-hand-side of this equality is equal to
$\langle (I_{\bfe_j}\otimes
\tilde{V}_{\bm}^*)\tilde{V}^{(j)*}\tilde{V}^{(i)}(I_{\bfe_i}\otimes
\tilde{V}_{\bn})(\xi\otimes h),\eta\otimes g \rangle =
\langle \tilde{V}_{\bm+\bfe_j}^*\tilde{V}_{\bn+\bfe_i}(\xi \otimes
h),\eta \otimes g \rangle =\langle V(\bn+\bfe_j-\bm-\bfe_i)(\xi
\otimes h),\eta \otimes g \rangle$ where the last equality follows
from Lemma~\ref{isom}. Thus, what we need to prove is
\be\label{a}
\langle (I_{\bfe_j}\otimes
\tilde{V}^{(i)})(I_{\bfe_i}\otimes \tilde{V}^{(j)*})
(I_{\bfe_i}\otimes
\tilde{V}_{\bn})(\xi\otimes h),(I_{\bfe_j}\otimes
\tilde{V}_{\bm})(\eta\otimes g) \rangle=
\ee
$$\langle V(\bn+\bfe_j-\bm-\bfe_i)(\xi
\otimes h),\eta \otimes g \rangle.$$
If $\bfe_j \leq \bn$ then the left-hand-side of the equation is
equal to $\langle (I_{\bfe_j}\otimes
\tilde{V}^{(i)})(I_{\bfe_i+\bfe_j}\otimes \tilde{V}_{\bn-\bfe_j})
(\xi\otimes h),(I_{\bfe_j}\otimes
\tilde{V}_{\bm})(\eta\otimes g) \rangle=\langle (I_{\bfe_j}\otimes
V(\bn-\bfe_j +\bfe_i)(\xi \otimes h),(I_{\bfe_j}\otimes
\tilde{V}_{\bm})(\eta\otimes g) \rangle=\langle V(\bn+\bfe_j-\bm-\bfe_i)(\xi
\otimes h),\eta \otimes g \rangle$. Similar argument works in the
case where $\bfe_i\leq \bm$. We now assume that $\bfe_j \nleq \bn$
and $\bfe_i \nleq \bm$. Hence $n_j=0=m_i$.

We first claim that, for $\bn\in \bZp$ and $j\in \{1,\ldots,k\}$
 with $\bn\wedge\bfe_j=0$,
we have
 $$\tilde{V}^{(j)*}\tilde{V}_{\bn}|X(\bn)\otimes H=(I_{\bfe_j}\otimes
\tilde{V}_{\bn})(I_{\bn}\otimes \tilde{V}_{\bfe_j}^*)|X(\bn)\otimes H.$$
Note that the ranges of the operators in this equation lie in
$E_j\otimes K$. Using (\ref{min}), (which is a consequence of the
minimality) it suffices to show, for every $\bp\in \bZp$, $\xi \in
X(\bn)$, $\eta \in X(\bp+\bfe_j)$ and $h,g \in H$,
\be
\langle \tilde{V}^{(j)*}\tilde{V}_{\bn}(\xi \otimes
h),(I_{\bfe_j}\otimes \tilde{V}_{\bp})(\eta \otimes g)\rangle =
\ee
$$\langle (I_{\bfe_j}\otimes
\tilde{V}_{\bn})(I_{\bn}\otimes \tilde{V}_{\bfe_j}^*)(\xi \otimes
h),(I_{\bfe_j}\otimes \tilde{V}_{\bp})(\eta \otimes g)\rangle.$$
Now, write $L$ for the left hand side of this equation and compute
$$L=\langle (I_{\bfe_j}\otimes \tilde{V}_{\bp}^*)\tilde{V}^{(j)*}
\tilde{V}_{\bn}(\xi \otimes
h),\eta \otimes g\rangle =\langle \tilde{V}_{\bp+\bfe_j}^*
\tilde{V}_{\bn}(\xi \otimes
h),\eta \otimes g\rangle =$$
$$\langle (I_{\bn\wedge
(\bp+\bfe_j)}\otimes V(\bn-\bp -\bfe_j)(\xi \otimes h),\eta
\otimes g\rangle=$$
$$\langle (I_{\bn\wedge
(\bp+\bfe_j)}\otimes T(\bn-\bp -\bfe_j)(\xi \otimes h),\eta
\otimes g\rangle$$
where the third equality follows from Lemma~\ref{isom} and for the
last one we use the
 regularity of the
dilation. Note that $(\bn-\bp -\bfe_j)_+=(\bn-\bp)_+$,
$(\bn-\bp-\bfe_j)_-=(\bn-\bp)_-+\bfe_j$, $\bn=\bn \wedge
(\bp+\bfe_j)+(\bn -\bp -\bfe_j)_+$ and $\bp +\bfe_j=\bn
\wedge(\bp+\bfe_j) + (\bn-\bp-\bfe_j)_-$.
 Thus, using
Lemma~\ref{consdc}(i), we have
 $T(\bn-\bp-\bfe_j)=(I_{(\bn-\bp-\bfe_j)_-}\otimes
 \tilde{T}_{(\bn-\bp-\bfe_j)_+})(I_{(\bn-\bp-\bfe_j)_+}\otimes
 \tilde{T}_{(\bn-\bp-\bfe_j)_-})$ and,
$$L=\langle (I_{\bp+\bfe_j}\otimes
\tilde{T}_{(\bn-\bp)_+})(I_{\bn}\otimes
\tilde{T}^*_{(\bn-\bp)_-+\bfe_j})(\xi \otimes h),\eta\otimes g
\rangle=$$
$$\langle (I_{\bp+\bfe_j}\otimes
\tilde{T}_{(\bn-\bp)_+})(I_{\bn+\bfe_j}\otimes
\tilde{T}^*_{(\bn-\bp)_-})(I_{\bn}\otimes \tilde{T}^*_{\bfe_j})
(\xi \otimes h),\eta\otimes g
\rangle.$$
Using Lemma~\ref{consdc}(i) again, this is equal to
$$\langle (I_{\bfe_j+\bn\wedge \bp}\otimes T(\bn-\bp))
(I_{\bn}\otimes \tilde{T}^*_{\bfe_j})
(\xi \otimes h),\eta\otimes g \rangle.$$
Since $\tilde{V}^*_{\bfe_j}h=\tilde{T}^*_{\bfe_j}h \in
X(\bfe_j)\otimes H$, for $h\in H$, and the dilation is regular,
this is equal to
$$\langle (I_{\bfe_j+\bn\wedge \bp}\otimes V(\bn-\bp))
(I_{\bn}\otimes \tilde{V}^*_{\bfe_j})
(\xi \otimes h),\eta\otimes g \rangle.$$
Applying Lemma~\ref{isom}, we get
$$L=\langle (I_{\bfe_j} \otimes \tilde{V}^*_{\bp}\tilde{V}_{\bn})
(I_{\bn}\otimes \tilde{V}^*_{\bfe_j})
(\xi \otimes h),\eta\otimes g \rangle=$$
$$\langle(I_{\bfe_j}\otimes
\tilde{V}_{\bn})(I_{\bn}\otimes \tilde{V}^*_{\bfe_j})
(\xi \otimes h),(I_{\bfe_j}\otimes \tilde{V}_{\bp})(\eta\otimes g )\rangle$$
proving the claim.

Now we turn to prove Equation (\ref{a}). The left hand side of
that equation is
$$\langle (I_{\bfe_i}\otimes \tilde{V}^{(j)*})
(I_{\bfe_i}\otimes
\tilde{V}_{\bn})(\xi\otimes h),(I_{\bfe_j}\otimes
\tilde{V}^{(i)*})(I_{\bfe_j}\otimes
\tilde{V}_{\bm})(\eta\otimes g) \rangle=$$
$$\langle (I_{\bfe_i}\otimes
(\tilde{V}^{(j)*}\tilde{V}_{\bn}))(\xi \otimes h), (I_{\bfe_j}\otimes
(\tilde{V}^{(i)*}\tilde{V}_{\bm}))(\eta \otimes g) \rangle.$$
Applying the claim, this is equal to
$$\langle (I_{\bfe_i+\bfe_j} \otimes
\tilde{V}_{\bn})(I_{\bn+\bfe_i}\otimes \tilde{V}_{\bfe_j}^*)(\xi
\otimes h),(I_{\bfe_i+\bfe_j} \otimes
\tilde{V}_{\bm})(I_{\bm+\bfe_j}\otimes \tilde{V}_{\bfe_i}^*)(\eta
\otimes g)=$$
$$\langle (I_{\bfe_i+\bfe_j+\bn\wedge \bm}\otimes
V(\bn-\bm))(I_{\bn+\bfe_i}\otimes \tilde{T}_{\bfe_j}^*)(\xi\otimes
h),(I_{\bm+\bfe_j}\otimes \tilde{T}_{\bfe_i}^*)(\eta\otimes
g)\rangle.$$
By regularity, this is equal to
$$\langle(I_{\bm+\bfe_j}\otimes \tilde{T}_{\bfe_i})
(I_{\bfe_i+\bfe_j+\bn\wedge \bm}\otimes
T(\bn-\bm))(I_{\bn+\bfe_i}\otimes \tilde{T}_{\bfe_j}^*)(\xi\otimes
h),(\eta\otimes
g)\rangle$$ and, applying Lemma~\ref{consdc}(i), we get
$$\langle(I_{\bm+\bfe_j}\otimes \tilde{T}_{\bfe_i})
(I_{\bfe_i+\bfe_j+ \bm}\otimes
\tilde{T}_{(\bn-\bm)_+})(I_{\bfe_i+\bfe_j+ \bn}\otimes
\tilde{T}_{(\bn-\bm)_-}^*)(I_{\bn+\bfe_i}\otimes \tilde{T}_{\bfe_j}^*)(\xi\otimes
h),$$
$$(\eta\otimes
g)\rangle=
\langle (I_{\bm+\bfe_j}\otimes \tilde{T}_{(\bn-\bm)_++\bfe_i})
(I_{\bn+\bfe_i}\otimes \tilde{T}_{(\bn-\bm)_-+\bfe_j}^*)(\xi
\otimes h),\eta\otimes g \rangle=$$
$$\langle (I_{\bm+\bfe_j}\otimes \tilde{T}_{(\bn+\bfe_i-\bm-\bfe_j)_+})
(I_{\bn+\bfe_i}\otimes \tilde{T}_{(\bn+\bfe_i-\bm-\bfe_j)_-}^*)(\xi
\otimes h),\eta\otimes g \rangle.$$
Using Lemma~\ref{consdc}(i) and the regularity of the dilation, we
find that the last expression is equal to
$$\langle V(\bn+\bfe_i-\bm-\bfe_j)(\xi \otimes h),\eta \otimes g
\rangle$$
proving (\ref{a}).

\end{proof}

\begin{lemma}\label{dcnc}
An isometric representation $(\rho,\{V^{(i)}\})$ is doubly commuting
if and only if, for every $\bn,\bm \in \bZp$,
\be\label{nc}
\tilde{V}_{\bn}\tilde{V}_{\bn}^*\tilde{V}_{\bm}\tilde{V}_{\bm}^*=\tilde{V}_{\bn\vee
\bm}\tilde{V}_{\bn\vee \bm}^* .
\ee
\end{lemma}
\begin{proof} Assume that the representation is doubly commuting
and compute, using Lemma~\ref{isom},
for $\bn,\bm \in \bZp$,
$$\tilde{V}_{\bn}\tilde{V}_{\bn}^*\tilde{V}_{\bm}\tilde{V}_{\bm}^*=
\tilde{V}_{\bn}(I_{\bn-(\bm-\bn)_-}\otimes V(\bm
-\bn))\tilde{V}_{\bm}^*.$$
Since the representation is doubly commuting, this is equal to
$$\tilde{V}_{\bn}(I_{\bn}\otimes
\tilde{V}_{(\bm-\bn)_+})(I_{\bm}\otimes
\tilde{V}_{(\bm-\bn)_-}^*)\tilde{V}_{\bm}^*=
\tilde{V}_{\bn+(\bm-\bn)_+}\tilde{V}_{\bm+(\bm-\bn)_-}^*=$$
$$\tilde{V}_{\bn
\vee \bm}\tilde{V}_{\bn\vee\bm}^*$$
proving one direction.
For the other direction, assume that (\ref{nc}) holds and fix
$i\neq j$ in $\{1,\ldots,k\}$. Then
$$\tilde{V}_{\bfe_j}(I_{\bfe_j}\otimes
\tilde{V}_{\bfe_i}\tilde{V}_{\bfe_i}^*)\tilde{V}_{\bfe_j}^*=
\tilde{V}_{\bfe_i+\bfe_j}\tilde{V}_{\bfe_i+\bfe_j}^*=\tilde{V}_{\bfe_i}
\tilde{V}_{\bfe_i}^*\tilde{V}_{\bfe_j}\tilde{V}_{\bfe_j}^*.$$
Multiplying on the left by $\tilde{V}_{\bfe_i}^*$ and on the right
by $\tilde{V}_{\bfe_j}$ and using the fact that the representation
is isometric, we get
$\tilde{V}_{\bfe_i}^*\tilde{V}_{\bfe_j}(I_{\bfe_j}\otimes
\tilde{V}_{\bfe_i}\tilde{V}_{\bfe_i}^*)=\tilde{V}_{\bfe_i}^*
\tilde{V}_{\bfe_j}$. Since $\tilde{V}_{\bfe_j}(I_{\bfe_j}\otimes
\tilde{V}_{\bfe_i})=\tilde{V}_{\bfe_i+\bfe_j}=
\tilde{V}_{\bfe_i}(I_{\bfe_i}\otimes
\tilde{V}_{\bfe_j})$, we have
$$\tilde{V}_{\bfe_i}^*
\tilde{V}_{\bfe_j}=\tilde{V}_{\bfe_i}^*(\tilde{V}_{\bfe_j}(I_{\bfe_j}\otimes
\tilde{V}_{\bfe_i}))(I_{\bfe_j}\otimes\tilde{V}_{\bfe_i}^*)=
\tilde{V}_{\bfe_i}^*(\tilde{V}_{\bfe_i}(I_{\bfe_i}\otimes
\tilde{V}_{\bfe_j}))(I_{\bfe_j}\otimes\tilde{V}_{\bfe_i}^*)=$$
$$(I_{\bfe_i}\otimes
\tilde{V}_{\bfe_j})(I_{\bfe_j}\otimes\tilde{V}_{\bfe_i}^*)$$
proving that the representation is doubly commuting.

\end{proof}

\begin{remark}\label{Nc}
An isometric representation satisfying (\ref{nc}) is referred to
in the literature as a Nica-covariant representation (see \cite{N} or \cite{F}).
Thus, the lemma shows that being Nica-covariant is equivalent to
being an isometric doubly commuting representation.
\end{remark}

An important representation of $X$ is the Fock representation.
It is
defined as in \cite{F}. We write
$$\mathcal{F}(X) = \sum_{\bn \in \bZp} \oplus X(\bn) .$$
As mentioned in \cite{F}, this is a $C^*$-correspondence over $A$
with left action given by
$$\varphi_{\infty}(a)(\oplus x_{\bn})=(\oplus \varphi_{\bn}(a)x_{\bn}).$$
We can define a representation $L$ of $X$ on $\mathcal{F}(X)$ by
setting
\be\label{L}
 L(x)(\oplus x_{\bn}) = \oplus (x\otimes x_{\bn})\;,\;\; \oplus x_{\bn} \in
\mathcal{F}(X).
\ee

Note that, strictly speaking, this is not what we defined as a
representation above (since $\mathcal{F}(X)$ is not a Hilbert space)
 but we can ``fix" it by representing
$\mathcal{L}(\mathcal{F}(X))$ on a Hilbert space.

Let $\mathcal{T}_c(X)$ be the $C^*$-algebra generated by the
operators $\{L(x): x\in X\}$.

If $\pi$ is a faithful representation of $A$ on a Hilbert space
$H$, then $\mathcal{F}(X)\otimes_{\pi}H$ is a Hilbert space and
the map $T\mapsto T\otimes I_H$ is a faithful representation of
$\mathcal{L}(\mathcal{F}(X))$ on $\mathcal{F}(X)\otimes_{\pi}H$
called \emph{the induced representation}. Its restriction to
$\mathcal{T}_c(X)$ is  a faithful representation of $\mathcal{T}_c(X)$
denoted $Ind(\pi)$.

In \cite[Theorem 6.3]{F}, Fowler proved the following.

\begin{theorem}\label{cov}(\cite{F})
There is a $C^*$-algebra, denoted $\mathcal{T}_{cov}(X)$, and an
 isometric representation $i_{X}:X\rightarrow \mathcal{T}_{cov}(X)$ such that
$\mathcal{T}_{cov}(X)$ is generated by $i_X(X)$ and
$(\mathcal{T}_{cov}(X),i_X)$ is universal for Nica-covariant
isometric representations of $X$, in the sense that:
\begin{enumerate}
\item[(a)] There is a faithful representation $\theta$ of
$\mathcal{T}_{cov}(X)$ on a Hilbert space such that $\theta \circ
i_X$ is a Nica-covariant isometric representation of $X$ ; and
\item[(b)] for every Nica-covariant isometric representation $(\sigma,T)$
of $X$, there is a $C^*$-representation $T\times \sigma$ of
$\mathcal{T}_{cov}(X)$ such that $T=(T\times \sigma)\circ i_X$.
\end{enumerate}
Up to canonical isomorphism, $(\mathcal{T}_{cov}(X),i_X)$ is the
unique pair with this property.
\end{theorem}

The following definition can be found in \cite[Definition 5.7]{F}.
Recall that, for a Hilbert $C^*$-module $E$, $\mathcal{K}(E)$ is
the closed ideal in $\mathcal{L}(E)$ generated by the (adjointable)
operators $\xi \otimes \eta^*$, for $\xi,\eta \in E$, defined by
$(\xi\otimes \eta^*)\zeta=\xi \langle \eta,\zeta \rangle$, $\zeta
\in E$.

\begin{definition}\label{compal}
We say that $X$ is compactly aligned if, whenever $T\in
\mathcal{K}(X(\bn))$ and $S\in \mathcal{K}(X(\bm))$, we have
$$ (S\otimes I_{\bn\vee \bm-\bm})(T\otimes I_{\bn\vee \bm
-\bn})\in \mathcal{K}(X(\bn\vee\bm)).$$
\end{definition}

Clearly, if, for every $\bn \in \bZp$,
$\mathcal{K}(X(\bn))=\mathcal{L}(X(\bn))$ then $X$ is compactly
aligned.

The proof of the following result can be dug out of \cite{F}.

\begin{theorem}\label{tcov}
Suppose $X$ is compactly aligned and each $X(\bn)$ ($\bn \in
\bZp$) is essential (that is, $\varphi_{X(\bn)}(A)X(\bn)$ is dense
in $X(\bn)$) then the pair ($\mathcal{T}_c(X),L)$ is canonically isomorphic to
($\mathcal{T}_{cov}(X),i_X)$. Thus, $(\mathcal{T}_c(X),L)$ is universal
for Nica-covariant (equivalently, for doubly commuting) isometric
representations of X.
\end{theorem}
\begin{proof}
Here we just indicate how to read the proof from the results of
\cite{F}. There, the author constructs a $C^*$-algebra denoted
$B_P\times_{\tau,X} P$ that contains $\mathcal{T}_{cov}(X)$
(Theorem 6.3 there). Let $\pi$ be a faithful nondegenerate representation
 of $A$ on a Hilbert space $H$ and write $\Psi$ for $Ind(\pi)\circ
 L$. This is an isometric, Nica-covariant, representation of $X$
 on $\mathcal{F}(X)\otimes_{\pi}H$
 (See Lemma 5.3 of \cite{F}). It gives rise to a representation, denoted
 $L^{\Psi}\times \Psi$, of $B_P\times_{\tau,X} P$ on
 $\mathcal{F}(X)\otimes_{\pi}H$ whose restriction to
 $\mathcal{T}_{cov}(X)$ is the Nica-covariant representation that
 $Ind(\pi)\circ L$ induces on $\mathcal{T}_{cov}(X)$ (by its
 universal property). In \cite[Corollary 7.7]{F} it is shown that
 $L^{\Psi}\times \Psi$ is a faithful representation. It follows
 that $Ind(\pi)\circ L$ gives rise to a faithful representation of
 $\mathcal{T}_{cov}(X)$ on $\mathcal{F}(X)\otimes_{\pi}H$. Its
 image is equal to the image of $Ind(\pi)$ and, thus, composing it
 with $Ind(\pi)^{-1}$, we get a $^*$-isomorphism from
 $\mathcal{T}_{cov}(X)$ onto $\mathcal{T}_{c}(X)$ that carries $i_X$
 to $L$.
 \end{proof}

\begin{definition}\label{tensor}
The Banach algebra generated by $\{L(x): x\in X(\bn),\bn\in \bZp
\}$ will be called \emph{the concrete tensor algebra} of $X$ and
will be written $\mathcal{T}_{+,c}(X)$.
\end{definition}

In \cite{S} we defined the tensor algebra $\mathcal{T}_+(X)$, associated
with $X$, as an algebra satisfying a certain universal property
(for c.c. representations of $X$). When $k=1$, it coincides with
the concrete tensor algebra $\mathcal{T}_{+,c}(X)$. In general,
the concrete tensor algebra does not have that universal property.
Nevertheless, it satisfies the following.

\begin{corollary}\label{repn}
Let $X$ be a compactly aligned product system of essential
correspondences.
For every c.c. doubly commuting representation $(\sigma,\{T^{(i)}\})$
of $X$ on a Hilbert space $H$, there is a completely contractive
representation $T\times \sigma$ of $\mathcal{T}_{+,c}(X)$ on $H$
such that, for every $1\leq i \leq k$ and every $\xi \in
X(\bfe_i)$,
$$(T\times \sigma)(L(\xi))=T^{(i)}(\xi) .$$
\end{corollary}
\begin{proof}
Let $(\rho,\{V^{(i)}\})$ be the minimal regular isometric dilation of
$(\sigma,\{T^{(i)}\})$ (on, say, $K$). By Theorem~\ref{doublycomm}
this isometric representation is doubly commuting. We see in
Lemma~\ref{dcnc} that it is Nica covariant. It then follows from
Theorem~\ref{tcov} that there is a $C^*$-representation
$\pi$ of $\mathcal{T}_c(X)$ on $K$ such that $V=\pi \circ L$.
Thus, for every $1\leq i \leq k$ and every $\xi \in
X(\bfe_i)$,
$\pi(L(\xi))=V^{(i)}(\xi) $. Writing $T \times\sigma $ for $P_H
\pi(\cdot)|H$, we see that $T\times \sigma$ is a completely
contractive map of $\mathcal{T}_c(X)$ into $B(H)$. Since $K\ominus
H$ is invariant under $V^{(i)}$ for all $1\leq i \leq k$ (and,
thus, invariant for $\tilde{V}_{\bn}$ for all $\bn \in \bZp$) the
map $T \times \sigma$ is multiplicative on $\mathcal{T}_{+,c}(X)$
and defines a completely contractive representation.

\end{proof}
\end{section}

\begin{section}{Examples}\label{ex}
\subsection{The case $k=1$}
In this case, we have a single $C^*$-correspondence $E$ over the
$C^*$-algebra $A$ and $X(n)=E^{\otimes n}$, $n \in \mathbb{Z}_+$. The algebra
$\mathcal{T}_{+,c}(X)$ was denoted by $\mathcal{T}_{+}(E)$ in
\cite{MS98} and its representations were studied there. Of course,
in this case, every representation is doubly commuting.
It was shown in \cite[Theorem 3.3]{MS98} that every c.c.
representation has a (unique) minimal isometric dilation. In
\cite[Theorem 3.10]{MS98} it was shown that every c.c.
representation of $E$ gives rise to a (unique) completely
contractive representation of $\mathcal{T}_+(E)$. Thus,
Theorem~\ref{doublycomm} and Corollary~\ref{repn} generalize these
results of \cite{MS98}.

\subsection{The case $A=E_i=\mathbb{C}$}
Now set $A=\mathbb{C}$ and, for each $1\leq i \leq k$,
$E_i=\mathbb{C}$ (with the obvious correspondence structure). In
order to define the product system $X$ (over $\bZp$) we need to
specify, for every $1\leq i,j \leq k$, an isomorphism of
correspondences $t_{i,j}:E_i \otimes E_j\rightarrow E_j \otimes E_i$ (with
$t_{j,i}=t_{i,j}^{-1}$ and $t_{i,i}=id$). This amounts to fixing complex numbers
$\lambda_{i,j} $ with $|\lambda_{i,j}|=1$, $\lambda_{i,i}=1$ and
$\lambda_{j,i}=\lambda_{i,j}^{-1}$ and setting $t_{i,j}(a\otimes
b)=\lambda_{i,j}b\otimes a$. (Note that (\ref{assoct}) is
satisfied).

 So, suppose we fix these
numbers and this defines $X$. Using (\ref{commute}),
 a c.c. representation of $X$ is now
a $k$-tuple $(T^{(1)},T^{(2)}, \ldots, T^{(k)})$ of contractions
in $B(H)$ (for some Hilbert space $H$) that satisfy
\be\label{1c}
T^{(i)}T^{(j)}=\lambda_{i,j}T^{(j)}T^{(i)}
\ee  for all $i,j$.
It is easy to check that this representation is doubly commuting
if and only if
\be\label{1dc}
T^{(i)*}T^{(j)}=\overline{\lambda_{i,j}}T^{(j)}T^{(i)*}
\ee for all $i\neq j$.

The case where $\lambda_{i,j}=1$ for all $i,j$ was studied
extensively and Theorem~\ref{doublycomm} and Corollary~\ref{repn}
are well known in this case (see, for example, \cite[Chapter I, Section
9]{SzF},
 \cite{GS} and \cite{Ti}). The algebra $\mathcal{T}_{+,c}(X)$ in
 this case is isomorphic to $A(\mathbb{D}^k)$ and
 Corollary~\ref{repn} amounts to the validity of the von Neumann inequality (for
 doubly commuting $k$-tuples).

If some of the $\lambda_{i,j}$'s are different from $1$,
$\mathcal{T}_{+,c}(X)$ is a non commutative subalgebra of
$B(l_2(\bZp))$. It is the Banach algebra generated by the
isometries $\{S_i: 1\leq i \leq k\}$ where (writing $\delta_{\bn}$
for the function in $l_2(\bZp)$ that is $1$ on $\bn$ and $0$
elsewhere)
 \be\label{Si}
 S_i \delta_{\bn}=\lambda(\bn,i)\delta_{\bn +\bfe_i}
 \ee
where $\lambda(\bn,i)=\prod_{j<i} \lambda_{j,i}^{n_j}$. (Note that
the isomorphism of $E_i \otimes X(\bn)=E_i \otimes E_1^{n_1}
\otimes E_2^{n_2} \otimes \cdots \otimes E_k^{n_k}$ and $X(\bn +
\bfe_i)=E_1^{n_1}
\otimes E_2^{n_2} \otimes \cdots \otimes E_i^{n_i+1} \otimes \cdots
 \otimes E_k^{n_k}$ sends $1\otimes 1 \otimes \cdots \otimes 1$ to
 $\lambda(\bn,i)(1\otimes \cdots \otimes 1)$).

For $\lambda:=\{\lambda_{i,j}\}$ as above, we write
$\mathcal{T}_{+,c}(\lambda)$ for the algebra
$\mathcal{T}_{+,c}(X)$ associated with the product system $X$ defined
by $\lambda$ (generated by the operators $S_i$ defined in (\ref{Si})).

The following Corollary is immediate from Theorem~\ref{doublycomm} and
Corollary~\ref{repn}. Part (ii) can be
viewed as a generalized von Neumann inequality.

\begin{corollary}\label{vni}
Fix $\lambda=\{\lambda_{i,j}:
|\lambda_{i,j}|=1,\;\lambda_{j,i}=\lambda_{i,j}^{-1}, \lambda_{i,i}=1 \}$ and let
$T^{(1)},T^{(2)},
\ldots, T^{(k)}$ be contractions in $B(H)$ that
satisfy (\ref{1c}) and (\ref{1dc}) above. Then
\begin{enumerate}
\item[(i)] there are isometries $U_1,U_2, \ldots, U_k$ (in $B(K)$, for
some Hilbert space $K$) that
satisfy (\ref{1c}) and (\ref{1dc}) and form a regular dilation of
$T^{(1)},T^{(2)},
\ldots \\ T^{(k)}$; and
 \item[(ii)] there is
a completely contractive representation $\pi$ of the algebra
$\mathcal{T}_{+,c}(\lambda)$  such that $\pi(S_i)=T^{(i)}$ for all
$1\leq i \leq k$. (Where $S_i$ are the operators defined in
(\ref{Si})). Thus, for every non commutative polynomial $p$ of $k$
variables,
$$\norm{p(T^{(1)},\ldots,T^{(k)})}\leq \norm{p(S_1,\ldots,S_k)}.$$
\end{enumerate}
\end{corollary}

\vspace{4mm}

If $dim H=1$, (\ref{1c}) implies (\ref{1dc}) and we get the
following.

\begin{corollary}\label{char}
The characters of $\mathcal{T}_{+,c}(\lambda)$ (that is, the one
dimensional representations of the algebra) are in one-to-one
correspondence with the set $\{t=(t_1,t_2,\ldots,t_k)\in
\mathbb{C}^k :\;|t_i|\leq 1 \;\emph{for all}\; 1\leq i \leq k,\;
t_it_j=0 \;\emph{whenever} \; \lambda_{i,j}\neq 1 \}$.
\end{corollary}

Now take $k=2$ and write $P_i=I-S_iS_i^*$. Then we have the
following.

\begin{corollary}\label{uv}
Let $k=2$ and assume that $\lambda :=\lambda_{1,2}$ is not a root of
unity. Let $J$ be the ideal of the $C^*$-algebra
$\mathcal{T}_c(X)$ generated by $P_1$ and $P_2$. Then
$\mathcal{T}_c(X)/ J$ is isomorphic to the irrational rotation
$C^*$-algebra $A_{\theta}$ (with $e^{2\pi i \theta}=\lambda$).
\end{corollary}
\begin{proof}
Write $q$ for the quotient map. Since $S_1S_2=\lambda S_2S_1$, the
same relation holds for $q(S_1)$ and $q(S_2)$. But these are
unitary operators and, thus, generate $A_{\theta}$.
\end{proof}

\subsection{The case $E_i=_{\alpha_i}A$}

Now fix a set of $k$ commuting $^*$-automorphisms
$\alpha_i$, $1\leq i \leq k$, of
$A$.
We write $_{\alpha_i}A$ for the $C^*$-correspondence over $A$
defined as follows. As a space, it is $A$. The left and right
actions  are defined by $\varphi(a)c
b=\alpha_i(a)c b$ (for $a,b\in A$ and $c\in _{\alpha_i}A$)
 and the inner product is
$\langle c_1,c_2
\rangle=c_1^*c_2$. Now let $E_i$ be $_{\alpha_i}A$. Note that,
for automorphisms $\alpha,\beta$ of $A$, $_{\alpha}A\otimes
_{\beta}A\cong _{\beta \alpha}A$ (via $a\otimes b \mapsto
\beta(a)b$). Since we assumed that the automorphisms $\alpha_i$
and $\alpha_j$ commute, we can combine these isomorphisms to get
an isomorphism $t_{i,j}:_{\alpha_i}A\otimes _{\alpha_j}A
\rightarrow _{\alpha_j}A\otimes _{\alpha_i}A$.
In fact, $t_{i,j}$ can be written explicitely: $t_{i,j}(a\otimes
b)=\alpha_i^{-1}\alpha_j(a)\otimes b$. It is easy to check that
condition (\ref{assoct}) holds and, therefore,
 this defines a
product system $X$.

Suppose $(\sigma,\{T^{(i)}\})$ is a c.c. representation of $X$ on $H$ with
a nondegenerate representation $\sigma$ of $A$. fix $i$ and a
(positive, contractive) approximate unit $\{u_{\lambda}\}$ in $A$
and consider, for $b\in A$,
$T^{(i)}(u_{\lambda})\sigma(b)=T^{(i)}(u_{\lambda}b)$. Since the
operators on the right converge (in norm, to $T^{(i)}(b)$), the
net  $\{T^{(i)}(u_{\lambda})\}$ has a strong operator limit $T_i$.
Then $T_i$ is a contraction and, for $b\in A$,
$T^{(i)}(b)=T_i\sigma(b)$. For every $a,b\in A$,
$T_i\sigma(\alpha_i(b))\sigma(a)=T^{(i)}(\alpha_i(b)a)=T^{(i)}(\varphi(b) a)=
\sigma(b)T^{(i)}(a)=\sigma(b)T_i\sigma(a)$. Thus, for every $b\in
A$,
\be\label{covalph}
T_i\sigma(\alpha_i(b))=\sigma(b)T_i.
\ee
Now, consider the commutation relation (\ref{commute}). Apply the
left hand side to $a\otimes b \otimes h \in _{\alpha_i}A
\otimes_{\alpha_j}A \otimes H$ to get $\tilde{T}^{(i)}(a\otimes
T^{(j)}(b)h)=T^{(i)}(a)T^{(j)}(b)h=T_i\sigma(a)T_j\sigma(b)h=
T_iT_j\sigma(\alpha_j(a)b)h$. Applying the right hand side to the
same element, we get $\tilde{T}^{(j)}(\alpha_i^{-1}\alpha_j(a)
\otimes
T^{(i)}(b)h)=T^{(j)}(\alpha_i^{-1}\alpha_j(a))T^{(i)}(b)h=\\
T_j\sigma(\alpha_i^{-1}\alpha_j(a))T_i\sigma(b)h=
T_jT_i\sigma(\alpha_j(a)b)h$. Thus the commutation relation is
equivalent to $T_iT_j=T_jT_i$ for every $i,j$. It follows that
every representation of $X$ is given by a (non degenerate)
representation $\sigma$ of $A$ on $H$ and by a $k$-tuple of
commuting contractions in $B(H)$ satisfying (\ref{covalph}).

Now we claim that such a representation is doubly commuting if and
only if the $k$-tuple is doubly commuting; that is,
$T_iT_j^*=T_j^*T_i$ for every $i\neq j$. To see this, first note
that, for $h\in H$ and $1\leq i \leq k$, $\tilde{T}^{(i)*}h$ is
the limit of $u_{\lambda}\otimes T_i^*h$. Indeed, for
$a\otimes g\in E_i \otimes_{\sigma}g$, we have $\langle
u_{\lambda}\otimes T_i^*h,a\otimes g \rangle = \langle
T_i^*h,\sigma(u_{\lambda})\sigma(a)g \rangle$ and, taking the
limit, we get $\langle T_i^*h,\sigma(a)g\rangle=
\langle h,T^{(i)}(a)g\rangle = \langle h,\tilde{T}^{(i)}(a\otimes
g)\rangle=
\langle
\tilde{T}^{(i)*}h,a \otimes g \rangle$.
To prove the claim we now apply the left hand side of equation
(\ref{eqdc}) to $a\otimes h \in E_i\otimes_{\sigma} H$
to get
$$\tilde{T}^{(j)*}T_i\sigma(a)h=\lim u_{\lambda}\otimes
T_j^*T_i\sigma(a)h=\lim u_{\lambda}\otimes
T_j^*\sigma(\alpha_i^{-1}(a))T_ih=$$
$$\lim u_{\lambda}\otimes
\sigma(\alpha_j\alpha_i^{-1}(a))T_j^*T_ih=\alpha_j\alpha_i^{-1}(a)\otimes
T_j^*T_ih.$$
Applying the right hand side of the same equation to $a\otimes h$
we get,
$$(I_{\bfe_j}\otimes \tilde{T}^{(i)})(t_{i,j}\otimes
I_H)(I_{\bfe_i}\otimes \tilde{T}^{(j)*})(a\otimes h)=\lim
(I_{\bfe_j}\otimes \tilde{T}^{(i)})(t_{i,j}\otimes
I_H)(a \otimes u_{\lambda}\otimes T_j^*h)$$
$$=\lim (I_{\bfe_j}\otimes
\tilde{T}^{(i)})(\alpha_i^{-1}\alpha_j(a)\otimes u_{\lambda}\otimes
T_j^*h)=\lim \alpha_i^{-1}\alpha_j(a)\otimes T^{(i)}(u_{\lambda})
T_j^*h=$$
$$\alpha_i^{-1}\alpha_j(a)\otimes T_iT_j^*h.$$
It follows that the representation is doubly commuting if and only
if the associated $k$-tuple is doubly commuting, as claimed.

In order to apply Corollary~\ref{repn}, note that, although
$\mathcal{L}(X(\bn))\neq \mathcal{K}(X(\bn))$ whenever $A$ is non
unital, the product system $X$ is easily seen to be compactly
aligned.

Now, it follows from Corollary~\ref{repn} that, given a representation
$\sigma$ of $A$ on $H$ and a doubly commuting $k$-tuple of
contractions $(T_1, \ldots, T_k)$ satisfying (\ref{covalph}), there is
a completely contractive representation of $\mathcal{T}_{+,c}(X)$
on $H$ sending $L(a)$ to $\sigma(a)$, if $a\in A=X(\bz)$, and to
$T_i\sigma(a)$ if $a\in _{\alpha_i}A=X(\bfe_i)$.

In order to relate the algebra $\mathcal{T}_{+,c}(X)$ to the
analytic crossed product studied in \cite{LM}, we write
$\gamma_i=\alpha_i^{-1}$, $1\leq i \leq k$, and note that
$\gamma_1,\ldots,\gamma_k$ define an action $\gamma$ of $\bZ$ on
$A$. The \emph{analytic crossed product} algebra
$A\times_{\alpha}\bZp$ is a subalgebra of the $C^*$
crossed product $A\times_{\alpha}\bZ$. The $C^*$ crossed product
is defined as the completion of the algebra $\ell^1(\bZ,A)$ (with
product defined by convolution and the involution and
$C^*$-norm are the natural ones). The analytic crossed product is
then
the Banach subalgebra
generated by the functions $\delta_{\bn,a}$ (for $\bn\in \bZp$ and
$a\in A$) defined by $\delta_{\bn,a}(\bm)=a$ if $\bn=\bm$ and $0$
otherwise.

For every $a\in A$, define $\sigma(a)=\delta_{\bz,a}$ and, for
$b\in _{\alpha_i}A$, set $T^{(i)}(b)=\delta_{\bfe_i,\alpha_i^{-1}(b)}$ to get an
isometric, doubly commuting, representation of $X$. It follows
from Theorem~\ref{tcov} that it yields a $C^*$-representation
$\pi$ of $\mathcal{T}_c(X)$ into (in fact, onto)
$A\times_{\gamma}\bZ$. Restricting $\pi$ to $\mathcal{T}_{+,c}$, we
get a completely contractive homomorphism
$$\pi_0:\mathcal{T}_{+,c}(X)\rightarrow A\times_{\gamma}\bZp.$$
Now let $\tau$ be a faithful (nondegenerate) representation of $A$
on a Hilbert space $H$ and write $V=Ind(\tau)\circ L$. By
\cite[Lemma 5.3]{F}, this is an isometric, Nica-covariant (hence,
doubly commuting) representation of $X$ on
$\mathcal{F}(X)\otimes_{\tau}H$. Using the results of \cite{LM},
it induces a completely contractive representation of
$A\times_{\gamma}\bZ$ on $\mathcal{F}(X)\otimes_{\tau}H$.
Combining it with $Ind(\tau)^{-1}$, we get a completely
contractive homomorphism
$$ \rho:A\times_{\gamma}\bZ \rightarrow \mathcal{T}_{+,c}(X).$$
Since $\rho$ and $\pi_0$ are the inverse of each other, we
conclude

\begin{corollary}\label{crossed}
For the product system $X$ defined by
$\alpha_1,\ldots,\alpha_k$ as above, the concrete tensor
algebra is completely isometrically isomorphic to the analytic
crossed product $A\times_{\gamma}\bZp$ where $\gamma$ is the
action induced by $\{\gamma_i=\alpha_i^{-1}\}$.
\end{corollary}

\begin{remark} The reason we need to consider
$A\times_{\gamma}\bZp$ instead of $A\times_{\alpha}\bZp$ can be
seen by comparing our covariance condition (\ref{covalph}) with
the covariance relation (1.2) in \cite{LM}.
\end{remark}

Finally, note that, in the construction of $X$ associated with
$\alpha_1,\alpha_2,\ldots,\alpha_k$ as above, we could also add a
``twist" to the multiplication, either by complex numbers (as in
Subsection~\ref{ex}.2) or by a family of unitaries in the center
of $A$ (satisfying a certain ``cocycle" identity that derives from
(\ref{assoct})).

\subsection{The case $A=\mathbb{C}$}

Now assume that $A=\mathbb{C}$ and, thus, each $E_i$ (and each
$X(\bm)$) is a Hilbert space. The isomorphisms $t_{i,j}:E_i
\otimes E_j \rightarrow E_j\otimes E_i$ are given by unitary
operators (satisfying the associativity condition (\ref{assoct})).
For simplicity, we assume here that each $E_i$ is finite
dimensional and write $d_i$ for its dimension and $\bd$ for
$(d_1,\ldots,d_k)$. (Note that the product system is compactly
aligned even in the infinite dimensional case). Also, we fix an
orthonormal basis $\{e^{(i)}_l: 1\leq l \leq
d_i\}$ for $E_i$ .

Note that  the algebra
$\mathcal{A}_{\bd,\theta}$, studied in \cite[Section
4]{KP}, is the
algebra $\mathcal{T}_{+,c}(X)$ defined in Definition~\ref{tensor} if each  $t_{i,j}$ is
induced from a permutation
 $\theta_{i,j}$ on $\{1,\ldots,d_i\}\times \{1,\ldots,d_j\}$
 in the sense that
\be\label{tperm}
 t_{i,j}(e^{(i)}_l\otimes
e^{(j)}_m)=e^{(j)}_r\otimes e^{(i)}_s
\ee
 whenever
$\theta_{i,j}(l,m)=(s,r)$. (And we write $\theta$ for the family
$\{\theta_{i,j}\}$ of these permutations, noting that it is assumed to
satisfy
an ``associativity" condition that can be derived from
(\ref{assoct})). In \cite[Theorem 4.1]{KP}, the authors studied
the one-dimensional representations of the algebra
$\mathcal{A}_{\bd,\theta}$ (that is, its characters). It is shown
there  that every one dimensional representation of $X$ gives rise
to such a character (and vice versa).

 For general representations
(not necessarily one dimensional) we restrict ourselves to the
doubly commuting ones. In order to present the consequences of
Theorem~\ref{doublycomm} and Corollary~\ref{repn} to the product
system $X$ with $A=\mathbb{C}$, we need the following definitions.

It will be convenient to write $[m]$ ($1\leq m\in \mathbb{Z}$) for
the set $\{1,\ldots,m\}$.

\begin{definition}\label{udc}
\begin{enumerate}
\item[(i)] A \emph{row contraction of length} $n$ on $H$ is an $n$-tuple
$T=(T_1,\ldots,T_n)$ of operators in $B(H)$ satisfying
$\sum^n_{i=1} T_iT_i^*\leq I$. Such a row contraction is a
\emph{row isometry} provided each $T_i$ is an isometry.
\item[(ii)] Let $u=(u_{(i,j)(l,p)})_{(i,j),(l,p)\in [n]\times
[m]}$ be a unitary matrix (of size $nm\times nm$), and $T$ and $S$ be row
contractions of lengths $n$ and $m$, respectively, on $H$.
 We say that the
(ordered) pair $(T,S)$ $u$-\emph{doubly commutes} if, for all
$1\leq i \leq n$ and $1\leq j \leq m$,
\begin{enumerate}
\item[(a)] $T_iS_j=\sum_{(p,l)\in [n]\times [m]}
u_{(i,j)(p,l)}S_lT_p$, and
\item[(b)] $S_j^*T_i=\sum_{(p,l)\in [n]\times [m]}
u_{(i,l)(p,j)}T_pS_l^*$.
\end{enumerate}
\end{enumerate}
\end{definition}

 Note that, once an orthonormal basis $\{e^{(i)}_l: 1\leq l \leq
d_i\}$ is fixed for every $E_i$, a unitary matrix $u$ of size
$d_id_j\times d_id_j$ as in Definition~\ref{udc}(ii), defines an
isomorphism $t$ from $E_i \otimes E_j$ onto $E_j\otimes E_i$ by
\be\label{tu}
t(e^{(i)}_q\otimes e^{(j)}_m)=\sum_{(p,l)\in [d_i]\times [d_j]}
u_{(q,m)(p,l)}e^{(j)}_l\otimes e^{(i)}_p.
\ee

\begin{theorem}\label{dctuples}
Let $\{u^{(i,j)}: i,j \in [k] \}$ be a family of unitary matrices
that define (via (\ref{tu})) a family $\{t_{i,j}\}$ of isomorphisms
 satisfying (\ref{assoct}) and let $(T^{(1)},\ldots,\\T^{(k)})$  be a
 $k$-tuple of row contractions on $H$ such that, for every $i\neq
 j$, $(T^{(i)},T^{(j)})$ $u^{(i,j)}$-doubly commutes.
Then it has a simultaneous (regular) dilation to a
$k$-tuple $(V^{(1)},\ldots,V^{(k)})$ of row isometries such that,
 for every $i\neq
 j$, $(V^{(i)},V^{(j)})$ $u^{(i,j)}$-doubly commutes.
\end{theorem}
\begin{proof}
The theorem follows immediately from Theorem~\ref{doublycomm} once
it is observed that each $T^{(i)}$ defines a c.c. representation
of $E_i$, condition (a) of Definition~\ref{udc}(ii) amounts to
condition (\ref{commute}) (that is, to the fact that the $k$-tuple
defines a representation of $X$) and condition (b) of
Definition~\ref{udc}(ii) amounts to the assumption that the
representation is doubly commuting.
\end{proof}

\begin{remark}\label{diag}
It is easy to check that, if each matrix $u^{(i,j)}$ above is
diagonal, condition (\ref{assoct}) is always satisfied.
\end{remark}

Applying Corollary~\ref{repn} we get.

\begin{corollary}\label{ccrepn}
Every $k$-tuple as in Theorem~\ref{dctuples} defines a completely
contractive representation of $\mathcal{T}_{+,c}(X)$.
\end{corollary}

Specializing to the situation studied in \cite[Section 4]{KP}, we
get the following.

\begin{corollary}\label{repnkgraph}
Suppose $\theta=\{\theta_{i,j}\}$ is a family of permutations as
in \cite{KP} (defining a product system $X$ via (\ref{tperm})) and
 $(T^{(1)},\ldots,T^{(k)})$ is a
 $k$-tuple of row contractions on $H$ such that, for every $i\neq
 j$ in $[k]$ and every $(l,m)\in [d_i]\times [d_j]$,
\begin{enumerate}
\item[(a)] $T^{(i)}_lT^{(j)}_m=T^{(j)}_sT^{(i)}_r$ where
$(r,s)=\theta_{i,j}(l,m)$, and
\item[(b)] $T^{(j)*}_mT^{(i)}_l =
\sum_{(r,m)=\theta_{i,j}(l,s)}T^{(i)}_rT^{(j)*}_s.$
\end{enumerate}
Then there is a completely contractive representation $\pi$ of
$\mathcal{A}_{\bd,\theta}$ (=$\mathcal{T}_{+,c}(X))$
 on $H$ mapping each $L_{e^{(i)}_l}$ (in the
notation of \cite{KP}, which is $L(e^{(i)}_l)$ in the sense of
 (\ref{L})) to $T^{(i)}_l$.

\end{corollary}

\end{section}

\end{document}